\documentclass[12pt]{article}
\usepackage{amsbsy,amsfonts,amsmath,amsthm,amssymb,euscript,enumerate}
\usepackage{graphicx,color,centernot,framed}
\definecolor{shadecolor}{gray}{0.85}
\usepackage{url}
\usepackage{hyperref}
\usepackage{ulem}
\hypersetup{colorlinks,citecolor=blue,linkcolor=blue,urlcolor=black}
\usepackage{bbm}
\newcommand{\ind}[1]{\mathbbm{1}{\raisebox{-2pt}{$\scriptstyle \{#1\}$}}}
\newcommand{\indic}[1]{\mathbbm{1}{\raisebox{-2pt}{$\scriptstyle #1$}}}

\allowdisplaybreaks[1]

\theoremstyle{plain}
\numberwithin{equation}{section}
\theoremstyle{plain}
\newtheorem{thm}{Theorem}[section]

\newtheorem{lmm}[thm]{Lemma}
\newtheorem{prp}[thm]{Proposition}

\theoremstyle{definition}

\newtheorem*{rem*}{Remark}

\newcommand{\Acal}{\mathcal{A}}

\newcommand{\db}{\Longleftrightarrow}
\newcommand{\dc}{d_\mathrm{c}}
\newcommand{\diff}{\mathrm{d}}
\newcommand{\dpst}{\displaystyle}

\newcommand{\lbeq}[1]{\label{eq:#1}}

\newcommand{\N}{{\mathbb N}}
\newcommand{\nn}{\nonumber}

\newcommand{\pc}{p_\mathrm{c}}
\newcommand{\Proof}[1]{\paragraph{\it #1}}
\newcommand{\QED}{\hspace*{\fill}\rule{7pt}{7pt}\bigskip}
\newcommand{\R}{{\mathbb R}}
\newcommand{\refeq}[1]{(\ref{eq:#1})}

\newcommand{\sss}{\scriptscriptstyle}
\newcommand{\Tcal}{\mathcal{T}}

\newcommand{\vno}{\varnothing}

\newcommand{\Z}{\mathbb{Z}}

\newcommand{\Zd}{\Z^d}

\makeatletter
\newcommand{\xRightarrow}[2][]{\ext@arrow 0359\Rightarrowfill@{#1}{#2}}
\makeatother

\definecolor{dgreen}{rgb}{0,0.6,0}

\title{Spread-out limit of the critical points for lattice trees and lattice 
animals in dimensions $d>8$}
\author{
Noe Kawamoto\footnote{Graduate School of Science, Hokkaido University, 
Japan.},~~~
Akira Sakai\footnote{Faculty of Science, Hokkaido University,
Japan. \url{https://orcid.org/0000-0003-0943-7842}}
}

\begin{document}
\maketitle
\begin{abstract}
A spread-out lattice animal is a finite connected set of edges in 
$\{\{x,y\}\subset\Zd:0<\|x-y\|\le L\}$.  A lattice tree is a lattice 
animal with no loops. The best estimate on the critical point $\pc$ so far was achieved by Penrose \cite{p94}: 
$\pc=1/e+O(L^{-2d/7}\log L)$ for both models for all $d\ge1$.  In this paper, 
we show that $\pc=1/e+CL^{-d}+O(L^{-d-1})$ for all $d>8$, where the 
model-dependent constant $C$ has the random-walk representation
\begin{align*}
C_\mathrm{LT}=\sum_{n=2}^\infty\frac{n+1}{2e}U^{*n}(o),&&
C_\mathrm{LA}=C_\mathrm{LT}-\frac1{2e^2}\sum_{n=3}^\infty U^{*n}(o),
\end{align*}
where $U^{*n}$ is the $n$-fold convolution of the uniform distribution on the 
$d$-dimensional ball $\{x\in\R^d:\|x\|\le1\}$.  The proof is based on a novel 
use of the lace expansion for $\tau_p(x)$ and detailed analysis of the 1-point 
function at a certain value of $p$ that is designed to make the analysis 
extremely simple.
\end{abstract}

\tableofcontents

\section{Introduction and the main result}\label{Introduction}
Given an $L\in\N$, we consider spread-out lattice animals $A=(V_A,E_A)$, where 
the vertex set $V_A$ is a finite subset of $\Zd$ and any pair of vertices in 
$V_A$ are connected by a path of spread-out edges 
$E_A\subset\{\{x,y\}:0<\|x-y\|\le L\}$; $\|\cdot\|$ is an arbitrary fixed norm on $\mathbb R^d$. A lattice tree is a lattice animal with no loops.  Both models are 
statistical-mechanical models for branched polymers.

\begin{figure}[t]
\begin{center}
\includegraphics[scale=0.6]{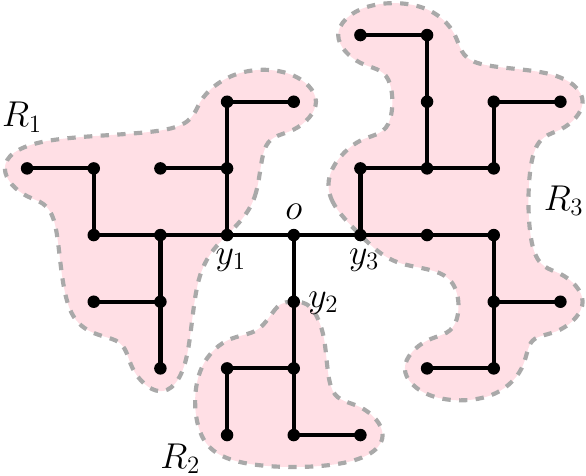}
\end{center}
\caption{A sample $T$ of $\Tcal_o$.  Removal of all edges $\{o,y_j\}\in E_T$ 
leaves disjoint subtrees $R_j$ rooted at $y_j$:
$V_T\setminus\{o\}=\bigcup_jV_{R_j}$ and 
$E_T\setminus\bigcup_j\{\{o,y_j\}\}=\bigcup_jE_{R_j}$.}
\label{fig:1pt}
\end{figure}

\begin{figure}[t]
\begin{center}
\includegraphics[scale=0.6]{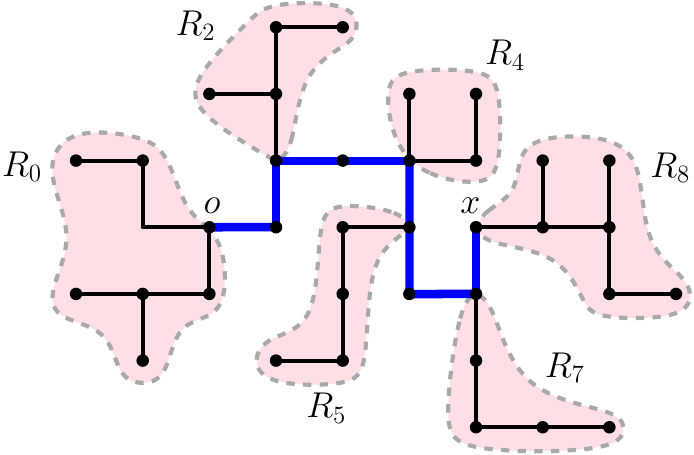}
\end{center}
\caption{A sample tree in $\Tcal_{o,x}$.  Removal of the backbone edges (in 
blue) yields disjoint subtrees $\{R_j\}$, called ribs.  In this example, 
$R_1,R_3$ and $R_6$ are single-vertex trees.}
\label{fig:2pt}
\end{figure}

To investigate their statistical properties, we consider the following 
generating functions.  Let 
\begin{align}\lbeq{Ddef}
\Lambda=\{x\in\Zd:0<\|x\|\le L\},&&
D(x)=\frac1{|\Lambda|}\ind{x\in\Lambda},
\end{align}
where $\indic{E}$ is the indicator function of $E$, being 1 or 0 depending on 
whether or not $E$ is true.  The function $D$ will be used as a transition 
probability of the underlying random walk.  Then, we define the weight function 
for a tree $T$ as
\begin{align}
W_p(T)=\prod_{\{x,y\}\in E_T}pD(x-y)=\bigg(\frac{p}{|\Lambda|}\bigg)^{|E_T|},
\end{align}
and similarly for a lattice animal $A$ as $W_p(A)$.  For a finite set 
$X\subset\Zd$, we denote by $\Tcal_X$ (resp., $\Acal_X$) the set of lattice 
trees $T$ with $X\subset V_T$ (resp., lattice animals $A$ with $X\subset V_A$); 
if $X$ consists of a vertex or two, we simply write, e.g., $\Tcal_o$ (for 
$X=\{o\}$; see Figure~\ref{fig:1pt}) 
or $\Tcal_{o,x}$ (for $X=\{o,x\}$; see Figure~\ref{fig:2pt}).  
The generating functions we want to 
investigate are the 1-point and 2-point functions, defiend respectively as
\begin{align}
g_p=\sum_{T\in\Tcal_o}W_p(T),&&
\tau_p(x)=\sum_{T\in\Tcal_{o,x}}W_p(T),
\end{align}
for lattice trees, and similarly defined for lattice animals.  
The susceptibility $\chi_p$ is the sum of the 2-point function, defined as 
\begin{align}\lbeq{chi-def}
\chi_p=\sum_{x\in\Zd}\tau_p(x)=\sum_{x\in\Zd}\sum_{T\in\Tcal_o}\ind{x\in V_T}\,
 W_p(T)=\sum_{T\in\Tcal_o}|V_T|\,W_p(T),
\end{align}
for lattice trees, and similarly for lattice animals.  It has been known (see 
\cite{ms93} and references therein) that there is a model-dependent critical 
point $\pc$ such that $\chi_p$ is finite if and only if $p<\pc$ and diverges as 
$p\uparrow\pc$.  The goal of this paper is to reveal the asymptotics of $\pc$ 
as $L\uparrow\infty$ for both models.

The best estimate so far on $\pc$ for the spread-out model was achieved by 
Penrose~\cite{p94}.  He investigated the growth constant, which is defined by 
the $n\uparrow\infty$ limit of the $n$th root of the number 
$t_n=\frac1n\sum_{T\in\Tcal_o}\ind{|V_T|=n}$ of $n$-vertex unrooted lattice 
trees.  Since $\{t_n\}_{n\in\N}$ is a supermultiplicative sequence, i.e., 
$t_{n+m}\ge t_n\,t_m$ (see, e.g., \cite{k81}), $\lim_{n\uparrow\infty}t_n^{1/n}$ 
exists and is asymptotically 
$e|\Lambda|+O(|\Lambda|^{5/7}\log|\Lambda|)$ as $|\Lambda|\uparrow\infty$ 
\cite{p94}.  Since $|E_T|=|V_T|-1$ for each lattice tree, we can rewrite 
$\chi_p$ as
\begin{align}
\chi_p\stackrel{\text{\refeq{chi-def}}}=\sum_{n=1}^\infty n\sum_{T\in\Tcal_o}
 \ind{|V_T|=n}\,W_p(T)=\sum_{n=1}^\infty n^2\bigg(\frac{p}{|\Lambda|}\bigg)^{n
 -1}t_n.
\end{align}
Therefore, for large $|\Lambda|$, 
\begin{align}\lbeq{Penrose}
\pc=\lim_{n\uparrow\infty}\bigg(\frac{n^2}{|\Lambda|^{n-1}}\,t_n\bigg)^{-1/n}
 =\frac{|\Lambda|}{e|\Lambda|+O(|\Lambda|^{5/7}\log|\Lambda|)}=\frac1e
 +O(|\Lambda|^{-2/7}\log|\Lambda|),
\end{align}
which is true for all dimensions $d\ge1$.  
Penrose also claimed in \cite[Section~3.1]{p94} that $\pc$ for lattice animals 
obeys the same bound, due to the result of Klarner~\cite{k67}.

A weaker estimate, $\pc=1/e+o(1)$ as $L\uparrow\infty$ for all $d$ bigger than 
the critical dimension $\dc=8$, was obtained by Miranda and Slade \cite{ms11}.  
In fact, their main concern was to obtain $1/d$ expansions of $\pc$ for the 
nearest-neighbor models.  In \cite{m12,ms13}, they showed that, 
\begin{align}\lbeq{nearest-neighbor}
\pc=\frac1e+\frac3{2e}|\Lambda|^{-1}+
 \begin{cases}
 \dfrac{115}{24e}|\Lambda|^{-2}+o(|\Lambda|^{-2})&[\text{lattice trees}],\\[1pc]
 \bigg(\dfrac{115}{24e}-\dfrac1{2e^2}\bigg)|\Lambda|^{-2}+o(|\Lambda|^{-2})
  &[\text{lattice animals}],
 \end{cases}
\end{align}
as $|\Lambda|=2d\uparrow\infty$. The proof is based on the lace expansion 
for the 2-point function $\tau_p(x)$ and an expansion for the 1-point function 
$g_p$ based on inclusion-exclusion. Notice that the model-dependence appears 
only from the $O(|\Lambda|^{-2})$ term, which is due to unit squares contained in $g_p$ for lattice animals, but not in $g_p$ for lattice trees. The lace expansion has been successful in 
showing mean-field critical behavior in high dimensions for various models, 
including lattice trees and lattice animals for $d>8$ (e.g., 
\cite{hhs03,hs90b,hs92b}).  The other models are self-avoiding walk for $d>4$ 
(e.g., \cite{bs85,hhs03,hs92a}), percolation for $d>6$ (e.g., 
\cite{hhs03,hs90a}), oriented percolation and the contact process for the 
spatial dimension $d>4$ (e.g., \cite{ny93,s01}), and the Ising and $\varphi^4$ 
models for $d>4$ (e.g., \cite{bhh21,s07,s15,s??}). 

For the nearest-neighbor lattice trees and lattice animals, in particular, Hara and Slade showed the mean-field critical behavior for both models in sufficiently high dimensions $d>8$ in \cite{hhs03,hs90b,hs92b}, where they did not mention about the specific dimension above which their results can be applied. In contrast, Fitzner and van der Hofstad proved in \cite{fh21} that the nearest-neighbor lattice trees and lattice animals obey the mean-field critical behavior above $d=16$ and $d=17$, respectively. The proof is based on the non-backtracking lace expansion (NoBLE) that is different from the standard lace expansion by Hara and Slade.

In \cite{hs05}, van der 
Hofstad and the second-named author of the current paper applied the lace 
expansion to the spread-out models (defined by $D$ in \refeq{Ddef}) of 
self-avoiding walk, percolation, oriented percolation and the contact process, 
and showed that, for all $d$ bigger than the respective critical dimension 
$\dc$,
\begin{align}
\pc=1+CL^{-d}+O(L^{-d-1}),
\end{align}
as $L\uparrow\infty$, where 1 is the mean-field value, and the model-dependent 
constant $C$ has the following random-walk representation:
\begin{align}
C=
 \begin{cases}
 \dpst\sum_{n=2}^\infty U^{*n}(o)
  &[\text{self-avoiding walk, the contact process}],\\[1pc]
 \dpst\frac12\sum_{n=2}^\infty U^{*2n}(o)&[\text{oriented percolation}],\\
 \dpst U^{*2}(o)+\sum_{n=3}^\infty\frac{n+1}2U^{*n}(o)&[\text{percolation}],
 \end{cases}
\end{align}
where $U^{*n}$ is the $n$-fold convolution in $\R^d$ of the uniform probability 
distribution $U$ on $\{x\in\R^d:\|x\|\le1\}$.  For example, if 
$\|x\|=\|x\|_\infty:=\max_j|x_j|$, then, for all $n\in\N$,
\begin{align}
U(x)=\frac{\ind{\|x\|_\infty\le1}}{2^d},&&
U^{*(n+1)}(x)=\int_{\R^d}U^{*n}(y)~U(x-y)~\mathrm{d}^dy.
\end{align}
These quantities are the spread-out limit of the underlying random walk 
generated by $D$ \cite[Section~4]{hs05}.  For example, for $d>4$,
\begin{align}\lbeq{RWs}
\sum_{n=2}^\infty\frac{n+1}2D^{*n}(o)=L^{-d}\sum_{n=2}^\infty\frac{n+1}2
 U^{*n}(o)+O(L^{-d-1}),
\end{align}
where we have used the same notation $*$ to represent convolutions on 
$\Zd$ as well. The error term $O(L^{-d-1})$ is due to Riemann-sum approximation.

We want to achieve a similar result for lattice trees and lattice animals, i.e., 
a random-walk representation for the difference between $\pc$ and its 
mean-field value $1/e$, and see how the model-dependence arises in it.  
In the rest of the paper, we will show the following:

\begin{shaded}
\begin{thm}\label{thm:main}
For both models with $d>8$ and $L\uparrow\infty$,
\begin{align}\lbeq{pc-asy}
\pc=\frac1e+CL^{-d}+O(L^{-d-1}),
\end{align}
where the model-dependent constant $C$ has the following random-walk 
representation:
\begin{align}\lbeq{coefficient}
C_\mathrm{LT}=\sum_{n=2}^\infty\frac{n+1}{2e}U^{*n}(o),&&
C_\mathrm{LA}=C_\mathrm{LT}-\frac1{2e^2}\sum_{n=3}^\infty U^{*n}(o),
\end{align}
\end{thm}
\end{shaded}
The difference in $\pc$ already shows up in the first error term of order $L^{-d}$ for the spread-out models, while it appears in \refeq{nearest-neighbor} from the second error term of order $d^{-2}$ for the nearest-neighbor models, as mentioned earlier. This is due to closed loops of length bigger than 2 in $g_p$ for lattice animals. The smallest among such loops for the spread-out model is of length 3 and of order $L^{-d}$, while that for the nearest-neighbor model is of length 4 and of order $d^{-2}$ (see Lemma~\ref{lmm:gp1a} below). Identifying coefficients of the higher order terms for the spread-out models may need more work since they are absorbed in the error term $O(L^{-d-1})$ in (1.12), which is inherent in Riemann-sum approximation, just as mentioned below \refeq{RWs}.

The proof of the above theorem is based on the lace expansion for the 2-point 
function and detailed analysis of the 1-point function, similarly to 
the previous work by Miranda and Slade \cite{ms13}. The different point of our method from that of them is to introduce a new base point $p_1$ defined in \refeq{p1gp1} below, as $p_1g_{p_1}=1$. It is to estimate various generating functions in terms of massless random walks. For the spread-out models of self-avoiding walk, percolation, oriented percolation and the contact process, van der Hofstad and Sakai \cite{hs05} simply used the base point $p_1=1$, because of the unity of the $1$-point function for those models. Since the analysis in terms of the underlying random walks is very simple, we do not have to know in detail the lace expansion; the exception is in Lemma~\ref{lmm:pc-2nd} below, where we investigate the first lace-expansion coefficient $\hat\pi_p^{\sss(1)}$ to prove $\pc-p_1=O(L^{-2d})$. However, the basic facts (summarized in Proposition~\ref{prp:xspace} below) and a {\it{minimum}} definition about the lace expansion coefficients should be enough to read the proof, which we hope makes this paper more accessible to wider audience. 

Our method can be applied to the nearest-neighbor model as well 
to identify the coefficient of $(2d)^{-1}$, as we can use the same method (i.e., Lemma~\ref{lmm:pc-2nd} below) to conclude $\pc-p_1=O(d^{-2})$, but this limit the accuracy our method can achieve. Therefore, to identify the higher-order coefficients, we may need investigate the lace expansion coefficients at $p_c$ more carefully as Miranda and Slade do in \cite{ms13}. 

The rest of the paper is organized as follows.  In Section~\ref{s:lace-exp}, 
we show that $\pc$ is close (up to order $L^{-2d}$) to $p_1$ that satisfies 
the identity $p_1g_{p_1}=1$, which is heavily used in the analysis in 
Sections~\ref{s:1pt} and \ref{s:difference}.  Section~\ref{s:1pt} is devoted to 
evaluating $g_{p_1}$ for lattice trees.  The 1-point function is split into two 
parts, $G$ and $H$, which are investigated in Sections~\ref{ss:G} and 
\ref{ss:H}, respectively.  Finally, in Section~\ref{s:difference}, we 
demonstrate how to evaluate the differnce between lattice trees and 
lattice animals.

\section{Results due to the lace expansion}\label{s:lace-exp}
In this section, we approximate $\pc$ by $p_1$ that is defined for both models 
by the identity
\begin{align}\lbeq{p1gp1}
p_1g_{p_1}=1.
\end{align}
From now on, we frequently use 
\begin{align}
\beta=L^{-d}.
\end{align}

\begin{shaded}
\begin{lmm}\label{lmm:pc-2nd}
For both models with $d>8$ and $L\uparrow\infty$,
\begin{align}\lbeq{pc-2nd}
0<\pc-p_1=O(\beta^2).
\end{align}
\end{lmm}
\end{shaded}

The key to the proof is the following collection of the lace-expansion results 
\cite{hhs03,l22}, in which we use
\begin{align}
h_p(x)=
 \begin{cases}
 0&[\text{lattice trees}],\\[5pt]
 \dpst(1-\delta_{o,x})\sum_{A\in\Acal_o}\ind{o\db x}\,W_p(A)\quad&[\text{lattice animals}],
 \end{cases}
\end{align}
where $o\db x$ means that ($o=x$ or) there is at least one pair of edge-disjoint 
paths from $o$ to $x$ in an animal $A$.  Let $\hat h_p$ denote the sum of 
$h_p(x)$ over $x\in\Zd$:
\begin{align}
\hat h_p=\sum_{x\in\Zd}h_p(x).
\end{align}

\begin{shaded}
\begin{prp}[\cite{hhs03,l22}]\label{prp:xspace}
For both models in $d>8$, there is a model-dependent $L_0<\infty$ such that, 
for all $L\ge L_0$, the following holds for all $p\le\pc$: 
\begin{enumerate}
\item The $1$-point function is bounded away from zero and infinity. In fact,
\begin{align}\lbeq{gp1O1}
1\le g_p\le4.
\end{align}
\item
There are nonnegative lace-expansion coefficients $\pi_p^{\sss(n)}(x)$, 
$n\in\N$, such that
\begin{align}\lbeq{pin-xbd}
\exists K<\infty,\quad\forall x\in\Zd,\quad
\pi_p^{\sss(n)}(x)\le\frac{KL^{-6}(K\beta)^{n-1}}{(\|x\|\vee L)^{2d-6}},
\end{align}
and that, by defining 
$\pi_p(x)=\sum_{n\in\N}(-1)^n\pi_p^{\sss(n)}(x)$, the recursion equation 
\begin{align}\lbeq{lace-exp}
\tau_p(x)&=g_p\delta_{o,x}+h_p(x)+\pi_p(x)\nn\\[5pt]
&\quad+\sum_{u,v}\Big(g_p\delta_{o,u}+h_p(u)+\pi_p(u)\Big)\,pD(v-u)\,
 \tau_p(x-v)
\end{align}
holds for all $x\in\Zd$.
\end{enumerate}
Consequently, there is a $K'<\infty$ such that
\begin{align}\lbeq{MFbehavior}
\forall x\ne o,\quad\tau_{\pc}(x)\le\frac{K'L^{-2}}{(\|x\|\vee L)^{d-2}},&&
\chi_p\underset{p\uparrow\pc}\asymp(\pc-p)^{-1/2}, 
\end{align}
where the latter means $\chi_p/(\pc-p)^{-1/2}$ is bounded 
away from 0 and $\infty$ as $p\uparrow\pc$, and 
\begin{align}\lbeq{pc-id}
\pc=\frac1{g_{\pc}+\hat h_{\pc}+\hat\pi_{\pc}}=\bigg(g_{\pc}+\sum_{x\ne o}
 h_{\pc}(x)+\sum_x\pi_{\pc}(x)\bigg)^{-1}.
\end{align}
\end{prp}
\end{shaded}

The above results for lattice trees are proven in \cite{l22} by following the 
same line of proof as in \cite{hhs03} and using the convolution bounds in 
\cite[Lemma~3.2]{cs15} instead of the weaker ones in 
\cite[Proposition~1.7]{hhs03}.  The same strategy applies to lattice animals, 
and we refrain from showing details.

Consequently, for any $p\le\pc$,
\begin{align}\lbeq{pin-hatbd}
\hat\pi_p^{\sss(n)}=\sum_{x\in\Zd}\pi_p^{\sss(n)}(x)
 \stackrel{\text{\refeq{pin-xbd}}}\le K(K\beta)^{n-1}\bigg(\sum_{x:\|x\|\le L}
 L^{-2d}+\sum_{x:\|x\|>L}\frac{L^{-6}}{\|x\|^{2d-6}}\bigg)=O(\beta)^n.
\end{align}
Moreover, by subadditivity (i.e., forgetting edge-disjointness among paths from 
$o$ to $x$),
\begin{align}\lbeq{hathbd}
\hat h_p\le\sum_{x\ne o}\tau_p(x)^2\stackrel{\text{\refeq{MFbehavior}}}\le
 (K'L^{-2})^2\bigg(\sum_{x:\|x\|\le L}L^{2(2-d)}+\sum_{x:\|x\|>L}\|x\|^{2(2-d)}
 \bigg)=O(\beta).
\end{align}
The identity \refeq{pc-id} is obtained by summing \refeq{lace-exp} over 
$x\in\Zd$, solving the resulting equation for $\chi_p$ and then using the fact 
that $\chi_p$ diverges as $p\uparrow\pc$. Substituting \refeq{pin-hatbd}-\refeq{hathbd} to \refeq{pc-id} yields\footnote{In 
\cite{l22}, Liang investigated $\hat\pi_{\pc}^{\sss(1)}$ in \refeq{pc-1st} for 
lattice trees and showed that, for all $d>8$,\; $\pc g_{\pc}$ rather than 
$\pc$ exhibits 
\begin{align}\lbeq{pcgpc}
\pc g_{\pc}=1+\frac\beta{e}\sum_{n=2}^\infty\binom{n+1}2U^{*n}(o)+O(\beta/L)
 \qquad\text{as }L\uparrow\infty.
\end{align}
This may be a bit of surprise, as the coefficient of $\beta$ is much larger 
than that in \refeq{pc-asy}--\refeq{coefficient}.
}
\begin{align}\lbeq{pc-1st}
\pc=\frac1{g_{\pc}}\bigg(1+\frac{\hat h_{\pc}-\hat\pi_{\pc}^{\sss(1)}}{g_{\pc}}
 +O(\beta^2)\bigg)^{-1}=\frac1{g_{\pc}}\bigg(1+\frac{\hat\pi_{\pc}^{\sss(1)}
 -\hat h_{\pc}}{g_{\pc}}\bigg)+O(\beta^2),
\end{align}
which is the starting point of the analysis.

\Proof{Proof of Lemma~\ref{lmm:pc-2nd}.}
First we show $p_1<\pc$.  Since $pg_p$ is increasing in $p$ with $p_1g_{p_1}=1$, 
it suffices to show $\pc g_{\pc}>1$. By \refeq{gp1O1} and \refeq{pc-1st}, it then suffices to 
show that $\hat\pi_{\pc}^{\sss(1)}-\hat h_{\pc}$ is bounded from below by 
$\beta$ times a positive constant for large $L$.  Here, and only here, we use 
the actual definition of the lace-expansion coefficient $\hat\pi_p^{\sss(1)}$ 
(see, e.g., \cite{hs90b}). We can easily check that $\hat\pi_p^{\sss(1)}$ for both 
models is larger than the sum of triangles consisting only of three 
distinct edges: $\hat\pi_p^{\sss(1)}\ge|\Lambda|(|\Lambda|-1)(p/|\Lambda|)^3$, 
which is enough for lattice trees because $\hat h _p \equiv 0$. For lattice animals, we show below $\hat h_p\le\frac14\hat\pi_p^{\sss(1)}+O(\beta^2)$ for $p \le p_c$ in high dimensions $d>8$. The aforementioned sufficient condition for $\pc g_{\pc}>1$  is now verified.  

Next we show $\pc-p_1=O(\beta^2)$ for lattice animals by induction. The same induction also works for lattice trees with $A=T$ and $\hat h_{p_c} \equiv 0$. Let $\{\ell_n\}_{n\in\N}$ be 
the following increasing sequence bounded above by 2:
\begin{align}\lbeq{hjdef}
\ell_1=1,&&
\ell_{j+1}=1+\frac{\ell_j}2\qquad[j\in\N].
\end{align}
Since $p_c=O(1)$ (see,e.g., \refeq{Penrose} or \cite[Proposition~2.2]{hhs03}) and $p_1g_{p_1}=\pc(g_{\pc}+\hat h_{\pc}+\hat\pi_{\pc})=1$, we have
\begin{align}\lbeq{1-p1/pc}
0<1-\frac{p_1}{\pc}=1-\frac{g_{\pc}+\hat h_{\pc}+\hat\pi_{\pc}}{g_{p_1}}
 =-\underbrace{\frac{g_{\pc}-g_{p_1}}{g_{p_1}}}_{\ge0}-\frac{\hat
 h_{\pc}+\hat\pi_{\pc}}{g_{p_1}},
\end{align}
which is bounded above by $-\hat\pi_{\pc}/g_{p_1}=O(\beta)$ (due to 
\refeq{pin-hatbd}), confirming $\pc-p_1=O(\beta^{\ell_1})$.

Now we suppose $\pc-p_1=O(\beta^{\ell_j})$. Notice that $g_{p_c}-g_{p_1}$ can be rewritten as
\begin{align}\lbeq{Fdef}
g_{\pc}-g_{p_1}=\sum_{A\in\Acal_o}\bigg(1-\Big(\frac{p_1}{\pc}\Big)^{|E_A|}
 \bigg)W_{\pc}(A)&=\bigg(1-\frac{p_1}{\pc}\bigg)\underbrace{\sum_{A\in\Acal_o}
 \sum_{n=0}^{|E_A|-1}\Big(\frac{p_1}{\pc}\Big)^nW_{\pc}(A)}_{=:F}\nn\\
&\!\!\!\stackrel{\text{\refeq{1-p1/pc}}}=\bigg(-\frac{g_{\pc}-g_{p_1}}{g_{p_1}}
 -\frac{\hat h_{\pc}+\hat\pi_{\pc}}{g_{p_1}}\bigg)F.
\end{align}
Solving this for $g_{\pc}-g_{p_1}$ yields
\begin{align}\lbeq{gc-g1}
g_{\pc}-g_{p_1}=-\frac{\hat h_{\pc}+\hat\pi_{\pc}}{g_{p_1}+F}F,
\end{align}
which is bounded above by $-\hat\pi_{\pc}=O(\beta)$ (due to \refeq{pin-hatbd}) 
for both models.  
By substituting \refeq{gc-g1} to \refeq{1-p1/pc}, we obtain
\begin{align}\lbeq{cancellation}
p_c-p_1=p_c\left(\frac{1}{g_{p_1}}\frac{\hat h_{\pc}+\hat\pi_{\pc}}{g_{p_1}+F}F-\frac{\hat
 h_{\pc}+\hat\pi_{\pc}}{g_{p_1}}\right)&=-p_c\frac{\hat h_{\pc}+\hat\pi_{\pc}}{g_{p_1}+F}\nn\\
&\overset{\refeq{pin-hatbd}}{=}p_c\frac{\hat\pi_{\pc}^{\sss (1)}-\hat h_{\pc}}{g_{p_1}+F}+O(\beta^2).
\end{align}
Recall the definition of $F$ in \refeq{Fdef}.  Since 
$(p_1/\pc)^nW_{\pc}(A)=(\pc/p_1)^{|E_A|-n}W_{p_1}(A)$, which is also true for 
lattice trees, we have
\begin{align}\lbeq{Flbd}
F=\sum_{A\in\Acal_o}\sum_{n=1}^{|E_A|}\Big(\frac{\pc}{p_1}\Big)^nW_{p_1}(A)
 \stackrel{p_1<\pc}\ge\sum_{A\in\Acal_o}|E_A|W_{p_1}(A)\stackrel{\sss|V_A|
 \le2|E_A|}\ge\frac{\chi_{p_1}}2.
\end{align}
By \refeq{cancellation} and \refeq{Flbd}, we can estimate $p_c-p_1$ as
\begin{align}\lbeq{bdpc-p1}
p_c-p_1=p_c\frac{\hat\pi_{\pc}^{\sss(1)}-\hat h_{\pc}}{g_{p_1}F^{-1}+1}F^{-1}+O(\beta^2)&= O(\beta)\chi_p^{-1}+O(\beta^2)\nn\\
&= O(\beta)(p_c-p_1)^{\frac{1}{2}}+O(\beta^2),
\end{align}
where, for the last inequality, we use $\chi_{p_1}\asymp(\pc-p_1)^{-1/2}$ for both 
models in dimensions $d>8$. Applying the inductive hypothesis $p_c-p_1=O(\beta^{\ell_j})$ to \refeq{bdpc-p1}, we obtain $p_c-p_1=O(\beta^{\ell_{j+1}})$. Therefore the induction completed. Since 
$\lim_{j\uparrow\infty}\ell_j=2$, this proves $\pc-p_1=O(\beta^2)$, as required.
\QED

\begin{figure}[t]
\begin{align*}
\includegraphics[scale=0.4]{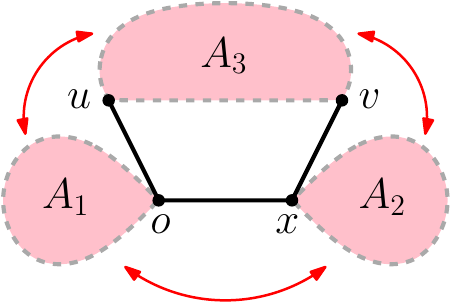}&&
\includegraphics[scale=0.4]{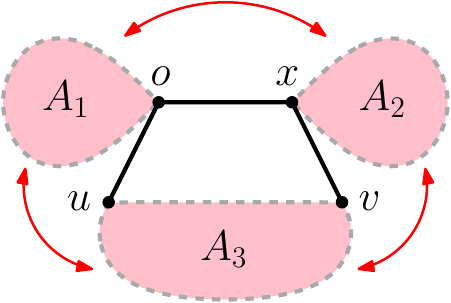}&&
\raisebox{-2pt}{\includegraphics[scale=0.4]{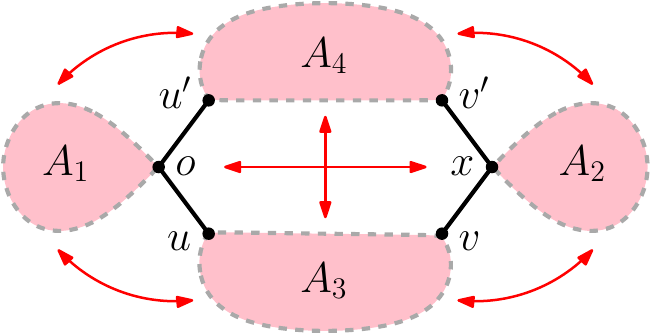}}
\end{align*}
\caption{Schematic representations of the three terms in \refeq{hdecomp}.  
The black line segments are pivotal for $o\db x$ in $A$.  Removal of those 
edges results in the animals $\{A_j\}_{j=1}^3$ or $\{A_j\}_{j=1}^4$ that are 
mutually avoiding, as indicated by the red arrows.  The vertices in $\Lambda$ 
are ordered in an arbitrary way (counter-clockwise in the above 2-dimensional 
figures).}
\label{fig:hdecomp}
\end{figure}

\Proof{Proof of $\hat h_p\le\frac14\hat\pi_p^{\sss(1)}+O(\beta^2)$ for lattice 
animals.}
First we recall that
\begin{align}
\hat h_p&=\sum_{x\ne o}\sum_{A\in\Acal_o}\bigg(\frac{p}{|\Lambda|}\bigg)^{|E_A|}
 \ind{o\db x}.
\end{align}
We split the sum into two depending on whether or not there are distinct 
vertices $y,z\in V_A$ such that 
$o\longleftrightarrow y$, $y\longleftrightarrow x$, $o\longleftrightarrow z$, 
$z\longleftrightarrow x$ and $y\longleftrightarrow z$ occur in $A$ 
edge-disjointly, i.e., those connections occur in distinct sets of $E_A$.  
(We note that, if $y=o$, for example, then we should interpret this as 
$o\longleftrightarrow x$, $o\db z$ and $z\longleftrightarrow x$ occuring 
edge-disjointly.)
Intuitively,
\begin{align}
\bigcup_{\substack{y,z\in V_A\\ (y\ne z)}}\raisebox{-1.5pc}{\includegraphics
 [scale=0.7]{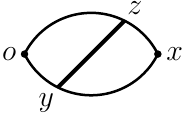}}.
\end{align}
Using submultiplicativity and the $x$-space bound in \refeq{MFbehavior}, 
we can show that the contribution from this case is $O(\beta^2)$.
On the other hand, if there are no such vertices $y,z\in V_A$, i.e., 
\begin{align}\lbeq{LAdoublyconn}
\raisebox{-1pc}{\includegraphics[scale=0.7]{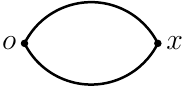}}~\setminus~
 \bigcup_{\substack{y,z\in V_A\\ (y\ne z)}}\raisebox{-1.5pc}{\includegraphics
 [scale=0.7]{LAdoublyconn2}},
\end{align}
then there are exactly two edge-disjoint connections between $o$ and $x$, 
with two pivotal edges from $o$, say $\{o,u\},\{o,u'\}$, and two from $x$, 
say $\{v,x\},\{v',x\}$, one of which may coincide with either $\{o,u\}$ or 
$\{o,u'\}$, for the double connection $o\db x$ in $A$.  Suppose that 
there is order among vertices in $\Lambda$.  If $u\in\Lambda$ is earlier than 
$u'\in\Lambda$ in this order, we write $u\prec u'$.  
Let $\Lambda(x)=\{v\in V:v-x\in\Lambda\}$.  Then, the contribution to 
$\hat h_p$ from \refeq{LAdoublyconn} is bounded above by
(see Figure~\ref{fig:hdecomp})
\begin{align}\lbeq{hdecomp}
\sum_{x\ne o}\sum_{\substack{A_1\in\Acal_o\\ A_2\in\Acal_x}}&W_p(A_1)\,W_p(A_2)
 \Bigg(\ind{x\in\Lambda}\frac{p}{|\Lambda|}\sum_{\substack{u\in\Lambda\\ (x\prec
 u)}}\sum_{\substack{v\in\Lambda(x)\\ (v\ne o)}}\bigg(\frac{p}{|\Lambda|}
 \bigg)^2\sum_{A_3\in\Acal_{u,v}}W_p(A_3)\nn\\
&+\ind{x\in\Lambda}\frac{p}{|\Lambda|}\sum_{\substack{u\in\Lambda\\ (u\prec x)}}
 \sum_{\substack{v\in\Lambda(x)\\ (v\ne o)}}\bigg(\frac{p}{|\Lambda|}\bigg)^2
 \sum_{A_3\in\Acal_{u,v}}W_p(A_3)\nn\\
&+\sum_{\substack{u,u'\in\Lambda\\ (u\prec u')}}\sum_{\substack{v,v'\in\Lambda
 (x)\\ (v\ne v')}}\bigg(\frac{p}{|\Lambda|}\bigg)^4\sum_{\substack{A_3\in
 \Acal_{u,v}\\ A_4\in\Acal_{u',v'}}}W_p(A_3)\,W_p(A_4)\Bigg)\prod_{i\ne j}
 \ind{V_{A_i}\cap V_{A_j}=\vno}.
\end{align}
Since $\Lambda$ is symmetric with respect to the underlying lattice symmetry, 
the first and second terms are the same.  Due to the same reason, the third 
term remains unchanged when the restriction $u\prec u'$ is replaced by 
$u'\prec u$.  Therefore, \refeq{hdecomp} equals
\begin{align}\lbeq{hdecomp-rewr}
\sum_{x\ne o}\sum_{\substack{A_1\in\Acal_o\\ A_2\in\Acal_x}}&W_p(A_1)\,W_p(A_2)
 \Bigg(\ind{x\in\Lambda}\frac{p}{|\Lambda|}\sum_{\substack{u\in\Lambda\\ (u\ne
 x)}}\sum_{\substack{v\in\Lambda(x)\\ (v\ne o)}}\bigg(\frac{p}{|\Lambda|}
 \bigg)^2\sum_{A_3\in\Acal_{u,v}}W_p(A_3)\nn\\
&+\frac12\sum_{\substack{u,u'\in\Lambda\\ (u\ne u')}}\sum_{\substack{v,v'\in
 \Lambda(x)\\ (v\ne v')}}\bigg(\frac{p}{|\Lambda|}\bigg)^4\sum_{\substack{A_3\in
 \Acal_{u,v}\\ A_4\in\Acal_{u',v'}}}W_p(A_3)\,W_p(A_4)\Bigg)\prod_{i\ne j}
 \ind{V_{A_i}\cap V_{A_j}=\vno}.
\end{align}

\begin{figure}[t]
\begin{center}
\includegraphics[scale=0.5]{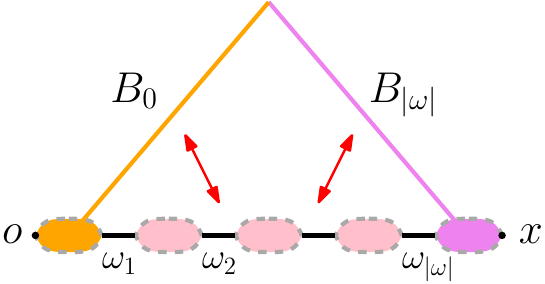}
\end{center}
\caption{Schematic representation of $\pi_p^{\sss(1)}(x)$.  
The sequence of edges $\omega_1,\dots,\omega_{|\omega|}$ joined by the animls 
$B_0,\dots,B_{|\omega|}$ form the backbone from $o$ to $x$ in $A$.  
The terminal animals $B_0$ and $B_{|\omega|}$ share a vertex (due to 
$\ind{B_0\cap B_{|\omega|}\ne\vno}$ in \refeq{pi1}), otherwise those animals 
are mutually avoiding (due to the product of indicators in \refeq{pi1}).  
Each animal $B_j$ contains a double connection between $\overline\omega_j$ 
and $\underline\omega_{j+1}$.}
\label{fig:pi1}
\end{figure}

Now we compare \refeq{hdecomp-rewr} with $\hat\pi_p^{\sss(1)}$ for lattice 
animals, which is defined as (see Figure~\ref{fig:pi1})
\begin{align}\lbeq{pi1}
\hat\pi_p^{\sss(1)}=\sum_x\sum_{\substack{\omega=\{(\underline\omega_i,
 \overline\omega_i)\}_{i=1}^{|\omega|}\\ (|\omega|\ge1)}}\bigg(\frac{p}
 {|\Lambda|}\bigg)^{|\omega|}\prod_{j=0}^{|\omega|}\sum_{B_j\in
 \Acal_{\overline\omega_j,\underline\omega_{j+1}}}W_p(B_j)\,\ind{\overline
 \omega_j\db\underline\omega_{j+1}\text{ in }B_j}\nn\\
\times\ind{B_0\cap B_{|\omega|}\ne\vno}\prod_{\substack{0\le k<l\le
 |\omega|\\ ((k,l)\ne(0,|\omega|))}}\ind{B_k\cap B_l=\vno},
\end{align}
where we have abused the notation $\overline\omega_0=o$ and 
$\underline\omega_{|\omega|+1}=x$.
This can be bounded below by restricting the sum over $\omega$ to those 
satisfying $\underline\omega_1=o$ and $\overline\omega_{|\omega|}=x$ 
(so that $\Acal_{\overline\omega_0,\underline\omega_1}=\Acal_o$ and 
$\Acal_{\overline\omega_{|\omega|},\,\underline\omega_{|\omega|+1}}=\Acal_x$) 
and then by restricting the sum over $B_0\in\Acal_o$ to $B_0=\{o\}$ (so that 
$\ind{B_0\cap B_{|\omega|}\ne\vno}=\ind{o\in B_{|\omega|}}$) or restricting 
the sum over $B_{|\omega|}\in\Acal_x$ to $B_{|\omega|}=\{x\}$ (so that 
$\ind{B_0\cap B_{|\omega|}\ne\vno}=\ind{x\in B_0}$):
\begin{align}\lbeq{pi1lowerevent}
\raisebox{-2pc}{\includegraphics[scale=0.45]{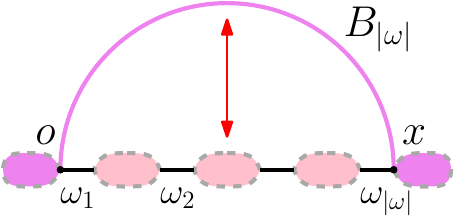}}\quad\cup\quad
 \raisebox{-2pc}{\includegraphics[scale=0.45]{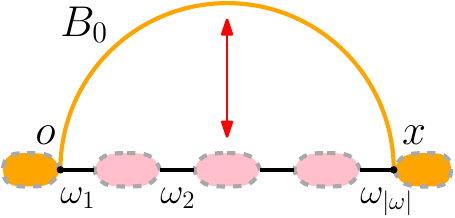}}.
\end{align}
Those two terms are basically the same.  Splitting the sum over $\omega$ into 
two depending on whether $|\omega|=1$ (so that $\omega=\{\omega_1\}$, where 
$\omega_1=(o,x)$) or $|\omega|\ge2$ and then, for 
the latter, by summing over the animals $B_1,\dots,B_{|\omega|-1}$ (to form 
an animal $A_3\in\Acal_{\overline\omega_1,\,\underline\omega_{|\omega|}}$), 
we obtain
\begin{align}\lbeq{pi1lb}
\hat\pi_p^{\sss(1)}\ge2\sum_{x\ne o}\sum_{B\in\Acal_{o,x}}&W_p(B)\Bigg(\ind{(o,
 x)\notin E_B}\ind{x\in\Lambda}\frac{p}{|\Lambda|}\nn\\
&+\sum_{\substack{u\in\Lambda\\ (u\ne x)}}\sum_{\substack{v\in\Lambda(x)\\
 (v\ne o)}}\bigg(\frac{p}{|\Lambda|}\bigg)^2\sum_{A_3\in\Acal_{u,v}}W_p(A_3)\,
 \ind{V_B\cap V_{A_3}=\vno}\Bigg).
\end{align}
We further bound this below by restricting the sum over $B\in\Acal_{o,x}$ to 
smaller animals $B=(V_B,E_B)$ with either
\begin{enumerate}[(i)]
\item
$V_B=V_{A_1}\cup V_{A_2}$, $E_B=E_{A_1}\cup\{(o,x)\}\cup E_{A_2}$ for some 
$A_1\in\Acal_o$, $A_2\in\Acal_x$ (as in the left and middle figures of 
Figure~\ref{fig:hdecomp}), or
\item
$V_B=V_{A_1}\cup V_{A_2}\cup V_{A_4}$, 
$E_B=E_{A_1}\cup\{(o,u')\}\cup E_{A_4}\cup\{(v',x)\}\cup E_{A_2}$ for some 
$A_1\in\Acal_o$, $A_2\in\Acal_x$, $u'\in\Lambda$, $v'\in\Lambda(x)$, 
$A_4\in\Acal_{u',v'}$ (as in the right figure of Figure~\ref{fig:hdecomp}).
\end{enumerate}
The contribution from (i) to the right-hand side of \refeq{pi1lb} is
\begin{align}\lbeq{pi1lbi}
2\sum_{x\ne o}\sum_{\substack{A_1\in\Acal_o\\ A_2\in\Acal_x}}W_p(A_1)\,W_p(A_2)
\ind{x\in\Lambda}\frac{p}{|\Lambda|}\sum_{\substack{u\in\Lambda\\ (u\ne x)}}
 \sum_{\substack{v\in\Lambda(x)\\ (v\ne o)}}\bigg(\frac{p}{|\Lambda|}\bigg)^2
 \sum_{A_3\in\Acal_{u,v}}W_p(A_3)\prod_{i\ne j}\ind{V_{A_i}\cap V_{A_j}=\vno},
\end{align}
while the contribution from (ii) is
\begin{align}\lbeq{pi1lbii}
2\sum_{x\ne o}\sum_{\substack{A_1\in\Acal_o\\ A_2\in\Acal_x}}&W_p(A_1)\,W_p(A_2)
 \Bigg(\ind{x\in\Lambda}\frac{p}{|\Lambda|}\sum_{\substack{u\in\Lambda\\ (u\ne
 x)}}\sum_{\substack{v\in\Lambda(x)\\ (v\ne o)}}\bigg(\frac{p}{|\Lambda|}
 \bigg)^2\sum_{A_3\in\Acal_{u,v}}W_p(A_3)\nn\\
&+\sum_{\substack{u,u'\in\Lambda\\ (u\ne u')}}\sum_{\substack{v,v'\in
 \Lambda(x)\\ (v\ne v')}}\bigg(\frac{p}{|\Lambda|}\bigg)^4\sum_{\substack{A_3\in
 \Acal_{u,v}\\ A_4\in\Acal_{u',v'}}}W_p(A_3)\,W_p(A_4)\Bigg)\prod_{i\ne j}
 \ind{V_{A_i}\cap V_{A_j}=\vno}.
\end{align}
Notice that the sum of \refeq{pi1lbi} and \refeq{pi1lbii} is four times as 
large as \refeq{hdecomp-rewr}.  This completes the proof of 
$\hat h_p\le\frac14\hat\pi_p^{\sss(1)}+O(\beta^2)$.
\QED

\section{Detailed analysis of the 1-point function for lattice trees}\label{s:1pt}
To complete the proof of Theorem~\ref{thm:main}, it remains to investigate 
$p_1=1/g_{p_1}$ (due to \refeq{p1gp1} and \refeq{pc-2nd}).  In this section, 
we concentrate our attention to lattice trees and show the following:

\begin{shaded}
\begin{lmm}\label{lmm:gp1}
For lattice trees with $d>8$ and $L\uparrow\infty$,
\begin{align}
g_{p_1}=e\bigg(1-\sum_{n=2}^\infty\frac{n+1}2D^{*n}(o)\bigg)
 +O(\beta^2).
\end{align}
Consequently,
\begin{align}
p_1=\frac1e+\sum_{n=2}^\infty\frac{n+1}{2e}D^{*n}(o)+O(\beta^2).
\end{align}
\end{lmm}
\end{shaded}

\begin{figure}[t]
\[ \raisebox{-3.6pc}{\includegraphics[scale=0.5]{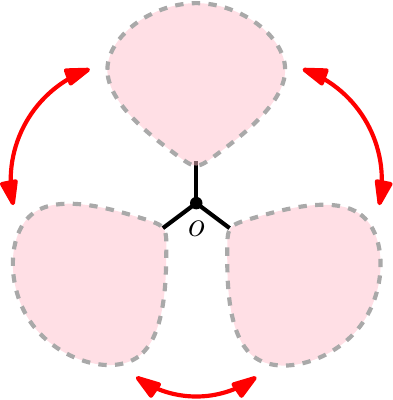}}~~~=~~~
 \raisebox{-3pc}{\includegraphics[scale=0.5]{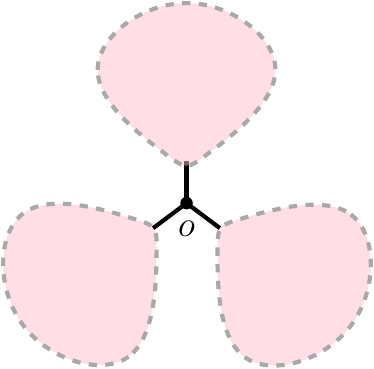}}~~~-~~~
 \raisebox{-3pc}{\includegraphics[scale=0.5]{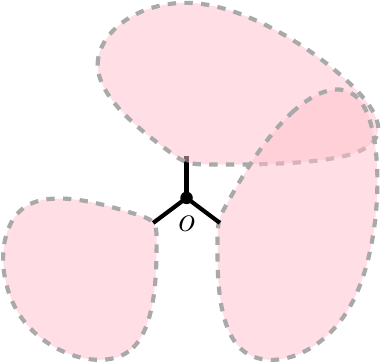}} \]
\caption{Intuitive explanation of \refeq{gp1=G-H}.  The double-headed arrows 
on the left ($=g_{p_1}$) represents mutual avoidance among subtrees.  In the 
first term on the right ($=G$), those subtrees are independently summed over 
$\Tcal_y\setminus\Tcal_o$, $y\in Y$, where $Y$ is the set of connected 
neighbors of the origin.  In the second term on the right ($=H$), there is at 
least one pair of subtrees that share vertices.}
\label{fig:gp1=G-H}
\end{figure}

To prove Lemma~\ref{lmm:gp1}, we first rewrite $g_{p_1}$ by identifying the 
connected neighbors $Y$ of the origin as
\begin{align}\lbeq{gp1-rewr}
g_{p_1}&=\sum_{T\in\Tcal_o}\bigg(\frac{p_1}{|\Lambda|}\bigg)^{|E_T|}\nn\\
&=1+\sum_{\substack{Y\subset\Lambda\\ (|Y|\ge1)}}\sum_{T\in\Tcal_{Y\cup\{o\}}}
 \bigg(\frac{p_1}{|\Lambda|}\bigg)^{|E_T|}\nn\\
&=1+\sum_{\substack{Y\subset\Lambda\\ (|Y|\ge1)}}\bigg(\frac{p_1}{|\Lambda|}
 \bigg)^{|Y|}\prod_{y\in Y}\sum_{R_y\in\Tcal_y\setminus\Tcal_o}\bigg(\frac{p_1}
 {|\Lambda|}\bigg)^{|E_{R_y}|}\prod_{\substack{u,v\in Y\\ (u\ne v)}}\ind{V_{R_u}
 \cap V_{R_v}=\vno},
\end{align}
where, and from now on, $\sum_{Y\subset\Lambda}$ is the sum over sets 
$Y$ of distinct vertices of $\Lambda$ (we recall that $o$ is not included in $\Lambda$) and 
\begin{align}\lbeq{a}
\prod_{y\in Y}\sum_{R_y\in\Tcal_y\setminus\Tcal_o}(p_1/|\Lambda|)^{|E_{R_y}|}
=\sum_{R_{y_1}\in\Tcal_{y_1}\setminus\Tcal_o}(p_1/|\Lambda|)^{|E_{R_{y_1}}|}
\cdots\sum_{R_{y_n}\in\Tcal_{y_n}\setminus\Tcal_o}(p_1/|\Lambda|)^{|E_{R_{y_n}}|}  
\end{align}
for $Y=\{y_1,\dots,y_n\}$. By convention, 
$\prod_{u\ne v}\ind{V_{R_u}\cap V_{R_v}=\vno}$ is regarded as 1 when $|Y|=1$.  
Let (see Figure~\ref{fig:gp1=G-H})
\begin{align}\lbeq{gp1=G-H}
g_{p_1}=G-H,
\end{align}
where
\begin{align}
G&=1+\sum_{\substack{Y\subset\Lambda\\ (|Y|\ge1)}}\bigg(\frac{p_1}{|\Lambda|}
 \bigg)^{|Y|}\prod_{y\in Y}\sum_{R_y\in\Tcal_y\setminus\Tcal_o}\bigg(\frac{p_1}
 {|\Lambda|}\bigg)^{|E_{R_y}|},\lbeq{Gdef}\\
H&=\sum_{\substack{Y\subset\Lambda\\ (|Y|\ge2)}}\bigg(\frac{p_1}{|\Lambda|}
 \bigg)^{|Y|}\prod_{y\in Y}\sum_{R_y\in\Tcal_y\setminus\Tcal_o}\bigg(\frac{p_1}
 {|\Lambda|}\bigg)^{|E_{R_y}|}\bigg(1-\prod_{\substack{u,v\in Y\\ (u\ne v)}}
 \ind{V_{R_u}\cap V_{R_v}=\vno}\bigg).\lbeq{Hdef}
\end{align}
We investigate those $G$ and $H$ in Sections~\ref{ss:G} and \ref{ss:H}, respectively 
(cf., Lemmas~\ref{lmm:G} and \ref{lmm:H} below).

\subsection{Detailed analysis of $G$}\label{ss:G}
From now on, we frequently use
\begin{align}\lbeq{Sdef}
S_{\ge t}(x)=\sum_{n=t}^\infty D^{*n}(x),
\end{align}
where $D^{*0}(x)=\delta_{o,x}$ by convention. The following is what we are going to show in this section:

\begin{shaded}
\begin{lmm}\label{lmm:G}
For lattice trees with $d>8$ and $L\uparrow\infty$,
\begin{align}\lbeq{Glmm}
G=e\bigg(1-\frac12D^{*2}(o)-S_{\ge2}(o)\bigg)+O(\beta^2).
\end{align}
\end{lmm}
\end{shaded}

\Proof{Proof.}
Since $p_1g_{p_1}=1$, we can rewrite $G$ as
\begin{align}
G&=1+\sum_{\substack{Y\subset\Lambda\\ (|Y|\ge1)}}\bigg(\frac{p_1}{|\Lambda|}
 \bigg)^{|Y|}\prod_{y\in Y}\Bigg(\underbrace{\sum_{R_y\in\Tcal_y}\bigg(
 \frac{p_1}{|\Lambda|}\bigg)^{|E_{R_y}|}}_{g_{p_1}}-\underbrace{\sum_{R_y\in
 \Tcal_{o,y}}\bigg(\frac{p_1}{|\Lambda|}\bigg)^{|E_{R_y}|}}_{\tau_{p_1}(y)}
 \Bigg)\nn\\
&=1+\sum_{\substack{Y\subset\Lambda\\ (|Y|\ge1)}}\bigg(\frac1{|\Lambda|}
 \bigg)^{|Y|}\prod_{y\in Y}\bigg(1-\frac{\tau_{p_1}(y)}{g_{p_1}}\bigg).
\end{align}
If we replace $\prod_{y\in Y}(1-\tau_{p_1}(y)/g_{p_1})$ by 1, 
then we obtain
\begin{align}\lbeq{G0-id}
G_0:=1+\sum_{\substack{Y\subset\Lambda\\ (|Y|\ge1)}}\bigg(\frac1{|\Lambda|}
 \bigg)^{|Y|}=\bigg(1+\frac1{|\Lambda|}\bigg)^{|\Lambda|}.
\end{align}
Since $k\log(1+1/k)=1-1/(2k)+O(k^{-2})$ as $k\uparrow\infty$, 
and since $|\Lambda|^{-1}=D^{*2}(o)$, we obtain
\begin{align}\lbeq{G0-est}
G_0=e\bigg(1-\frac1{2|\Lambda|}\bigg)+O(|\Lambda|^{-2})
 =e\bigg(1-\frac12D^{*2}(o)\bigg)+O(\beta^2).
\end{align}

Next we consider the remainder, which is
\begin{align}
G-G_0&=\sum_{\substack{Y\subset\Lambda\\ (|Y|\ge1)}}\bigg(\frac1{|\Lambda|}
 \bigg)^{|Y|}\Bigg(\prod_{y\in Y}\bigg(1-\frac{\tau_{p_1}(y)}{g_{p_1}}\bigg)-1
 \Bigg)\nn\\
&=\sum_{\substack{Y\subset\Lambda\\ (|Y|\ge1)}}\bigg(\frac1{|\Lambda|}
 \bigg)^{|Y|}\sum_{\substack{Z\subset Y\\ (|Z|\ge1)}}\prod_{y\in Z}\frac{-
 \tau_{p_1}(y)}{g_{p_1}}.
\end{align}
Changing the order of sums yields
\begin{align}\lbeq{G-G0=G1+G2}
G-G_0&=\sum_{\substack{Z\subset\Lambda\\ (|Z|\ge1)}}\prod_{y\in Z}\frac{-
 \tau_{p_1}(y)}{g_{p_1}}\bigg(\frac1{|\Lambda|}\bigg)^{|Z|}\sum_{Z\subset Y
 \subset\Lambda}\bigg(\frac1{|\Lambda|}\bigg)^{|Y\setminus Z|}\nn\\
&=\sum_{\substack{Z\subset\Lambda\\ (|Z|\ge1)}}\prod_{y\in Z}\frac{-\tau_{p_1}
 (y)\,D(y)}{g_{p_1}}\bigg(1+\frac1{|\Lambda|}\bigg)^{|\Lambda
 \setminus Z|}\nn\\
&=G_1+G_2,
\end{align}
where
\begin{align}
G_1&=\sum_{\substack{Z\subset\Lambda\\ (|Z|=1)}}\prod_{y\in Z}\frac{-\tau_{p_1}
 (y)\,D(y)}{g_{p_1}}\bigg(1+\frac1{|\Lambda|}\bigg)^{|\Lambda\setminus Z|}
 =\frac{-(\tau_{p_1}*D)(o)}{g_{p_1}}\frac{G_0}{(1+1/|\Lambda|)},\lbeq{G1def}\\
G_2&=\sum_{\substack{Z\subset\Lambda\\ (|Z|\ge2)}}\prod_{y\in Z}\frac{-
 \tau_{p_1}(y)\,D(y)}{g_{p_1}}\bigg(1+\frac1{|\Lambda|}\bigg)^{|\Lambda\setminus
 Z|}.\lbeq{G2def}
\end{align}
We make use of reflection symmetry for \refeq{G-G0=G1+G2}--\refeq{G2def}, and will use it below frequently. To estimate $G_1$ and $G_2$, we use the following lemma, 
which will be proven after the proof of Lemma~\ref{lmm:G} is completed.

\begin{shaded}
\begin{lmm}\label{lmm:RWapprox}
For any $d>2$ and $x\ne o$, the lattice-tree 2-point function satisfies
\begin{align}\lbeq{RWapprox}
0\le S_{\ge1}(x)-\frac{\tau_{p_1}(x)}{g_{p_1}}\le\sum_{\substack{y,z\in\Zd\\
 (y\ne z)}}S_{\ge0}^{*2}(z-y)\,S_{\ge0}(y)\,S_{\ge1}(z-y)\,S_{\ge0}(x-z).
\end{align}
\end{lmm}
\end{shaded}
\begin{rem*}
The right-hand side of \refeq{RWapprox} is diagrammatically represented by 
\begin{align}\lbeq{RWapproxdiagramA}
\sum_{\substack{y,z\in\Zd\\
 (y\ne z)}}S_{\ge0}^{*2}(z-y)\,S_{\ge0}(y)\,S_{\ge1}(z-y)\,S_{\ge0}(x-z)=~
 \raisebox{-0.5pc}{\includegraphics[scale=0.3]{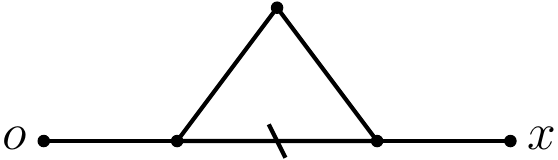}},
\end{align}
where an unslashed (resp., slashed) line represents $S_{\ge0}$ (resp., 
$S_{\ge1}$) and an unlabelled vertex is summed over $\Zd$. Due to translation invariance, we can change the order of 
terms in \refeq{RWapprox} for a given $x\ne o$.  Then \refeq{RWapproxdiagramA} is also equal to
\begin{align}\lbeq{RWapproxdiagramB}
\sum_{\substack{y\in\Zd\\ (y\ne x)}}S_{\ge0}^{*2}(y)\,S_{\ge0}^{*2}(x-y)\,
 S_{\ge1}(x-y)=~\raisebox{-0.5pc}{\includegraphics[scale=0.3]{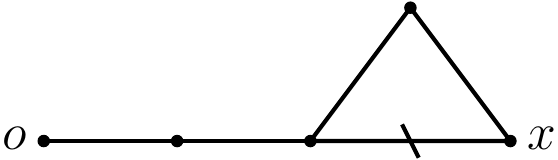}}.
\end{align}
These diagramatic representations will be used in the proof of Lemma~\ref{lmm:H} below.
\end{rem*}

First we estimate $G_2$.  By the first inequality in \refeq{RWapprox} and 
the heat-kernel bound (see, e.g., \cite[(1.6)]{cs08}):
\begin{align}\lbeq{D*bd}
\|D^{*n}\|_\infty=O(\beta)n^{-d/2}\qquad[n\in\N],
\end{align}
we can show
\begin{align}\lbeq{roughbd}
\frac{(\tau_{p_1}*D)(o)}{g_{p_1}}\le\sup_{x\in\Lambda}\frac{\tau_{p_1}(x)}
 {g_{p_1}}\le\sup_{x\in\Lambda}S_{\ge1}(x)\le\sum_{n=1}^\infty\|D^{*n}\|_\infty
 \stackrel{d>2}=O(\beta).
\end{align}
Therefore,
\begin{align}\lbeq{G2-est}
|G_2|&\le G_0\sum_{n=2}^{|\Lambda|}\sum_{\substack{Z\subset\Lambda\\ (|Z|=n)}}
 \prod_{y\in Z}\frac{\tau_{p_1}(y)\,D(y)}{g_{p_1}}\le G_0\sum_{n=2}^\infty
 \bigg(\frac{(\tau_{p_1}*D)(o)}{g_{p_1}}\bigg)^n=O(\beta^2).
\end{align}

Next we estimate $G_1$ in \refeq{G1def}.  By using \refeq{RWapprox}, 
we have

\begin{align}\lbeq{b}
0\le S_{\ge2}(o)-\frac{(\tau_{p_1}*D)(o)}{g_{p_1}}&=\sum_{x \in \Zd}D(x)\left(S_{\ge1}(x)-\frac{\tau_{p_1}(x)}{g_{p_1}}\right)\nn\\
&\overset{\refeq{RWapprox}}{\le} \sum_{x \in \Zd}D(x)\sum_{\substack{y,z\in\Zd\\
 (y\ne z)}}S_{\ge0}^{*2}(z-y)\,S_{\ge0}(y)\,S_{\ge1}(z-y)\,S_{\ge0}(x-z)\nn\\
&\le\sum_{\substack{y,z\in\Zd\\ (y\ne z)}}S_{\ge0}^{*2}(z-y)\,S_{\ge0}(y)\,
 S_{\ge1}(z-y)\,S_{\ge1}(z)\nn\\
&\le\sup_{w\ne o}S_{\ge0}^{*2}(w)~\big(S_{\ge0}*S_{\ge1}^{*2}\big)(o).
\end{align}
By the heat-kernel bound \refeq{D*bd}, we can estimate each term as
\begin{align}
S_{\ge0}^{*2}(w)&=\sum_{s,t=0}^\infty D^{*(s+t)}(w)\stackrel{w\ne o}\le
 \sum_{n=1}^\infty(n+1)D^{*n}(w)\stackrel{d>4}=O(\beta),\lbeq{RWanal1}\\
\big(S_{\ge0}*S_{\ge1}^{*2}\big)(o)&=\sum_{s=0}^\infty\sum_{t,u=1}^\infty
 D^{*(s+t+u)}(o)=\sum_{n=2}^\infty\binom{n}2D^{*n}(o)\stackrel{d>6}=O(\beta),
 \lbeq{RWanal2}
\end{align}
so that
\begin{align}\lbeq{RWest1LT}
0\le S_{\ge2}(o)-\frac{(\tau_{p_1}*D)(o)}{g_{p_1}}=O(\beta^2).
\end{align}
Therefore, by \refeq{G0-est},
\begin{align}\lbeq{G1-est}
G_1&=(-\underbrace{S_{\ge2}(o)}_{\overset{\refeq{roughbd}}{=}O(\beta)}+O(\beta^2))\bigg(e-\underbrace{\frac{e}{2}D^{*2}(o)}_{\overset{\refeq{D*bd}}{=}O(\beta)}+O(\beta^2)\bigg)\nn\\
&=-eS_{\ge2}(o)+O(\beta^2).
\end{align}
Summarizing \refeq{G0-est}, \refeq{G-G0=G1+G2}, \refeq{G2-est} and 
\refeq{G1-est}, we complete the proof of Lemma~\ref{lmm:G}.
\QED

\Proof{Proof of Lemma~\ref{lmm:RWapprox}.}
First we recall
\begin{align}
\frac{\tau_{p_1}(x)}{g_{p_1}}=\frac1{g_{p_1}}\sum_{T\in\Tcal_{o,x}}
 \bigg(\frac{p_1}{|\Lambda|}\bigg)^{|E_T|}.
\end{align}
Since a tree $T\in\Tcal_{o,x}$ can be divided into a unique path 
$\omega=(\omega_0,\omega_1,\dots,\omega_{|\omega|})$ from $\omega_0=o$ to 
$\omega_{|\omega|}=x$, called a backbone, and disjoint subtrees 
$R_j\in\Tcal_{\omega_j}$, called ribs (see Figure~\ref{fig:2pt}), 
we can rewrite the above expression as
\begin{align}\lbeq{RWapprox-ubd1}
\frac{\tau_{p_1}(x)}{g_{p_1}}=\frac1{g_{p_1}}\sum_{\omega:o\to x}\bigg(
 \frac{p_1}{|\Lambda|}\bigg)^{|\omega|}\prod_{j=0}^{|\omega|}\sum_{R_j\in
 \Tcal_{\omega_j}}\bigg(\frac{p_1}{|\Lambda|}\bigg)^{|E_{R_j}|}\prod_{s<t}
 \ind{V_{R_s}\cap V_{R_t}=\vno}.
\end{align}
If we replace the indicator $\prod_{s<t}\ind{V_{R_s}\cap V_{R_t}=\vno}$ by 1, 
then we obtain
\begin{align}\lbeq{RWapprox-ubd2}
\frac1{g_{p_1}}\sum_{\omega:o\to x}\bigg(\frac{p_1}{|\Lambda|}\bigg)^{|\omega|}
 \prod_{j=0}^{|\omega|}\underbrace{\sum_{R_j\in\Tcal_{\omega_j}}\bigg(\frac{p_1}
 {|\Lambda|}\bigg)^{|E_{R_j}|}}_{g_{p_1}}\stackrel{\text{\refeq{p1gp1}}}
 =\sum_{\omega:o\to x}\bigg(\frac1{|\Lambda|}\bigg)^{|\omega|}
 \stackrel{x\ne o}=S_{\ge1}(x).
\end{align}

Next we consider the remainder.  Since $1-\prod_{j=1}^na_j\le \sum_{j=1}^n(1-a_j)$ as long as $0\le a_j\le1$ for all $j$, we can bound the remainder as
\begin{align}\lbeq{c}
S_{\ge1}(x)-\frac{\tau_{p_1}(x)}{g_{p_1}}&=\frac1{g_{p_1}}
 \sum_{\omega:o\to x}\bigg(\frac{p_1}{|\Lambda|}\bigg)^{|\omega|}\prod_{j
 =0}^{|\omega|}\sum_{R_j\in\Tcal_{\omega_j}}\bigg(\frac{p_1}{|\Lambda|}
 \bigg)^{|E_{R_j}|}\bigg(1-\prod_{s<t}\ind{V_{R_s}\cap V_{R_t}=\vno}\bigg)\nn\\
&\le\frac1{g_{p_1}}\sum_{\omega:o\to x}\bigg(\frac{p_1}{|\Lambda|}
 \bigg)^{|\omega|}\prod_{j=0}^{|\omega|}\sum_{R_j\in\Tcal_{\omega_j}}\bigg(
 \frac{p_1}{|\Lambda|}\bigg)^{|E_{R_j}|}\sum_{s<t}\ind{V_{R_s}\cap V_{R_t}\ne
 \vno}.
\end{align}
If $V_{R_s}\cap V_{R_t}\ne\vno$, then there must be a $w\in\Zd$ that is 
shared by those two ribs.  Therefore, the remainder is further bounded above as
\begin{align}\lbeq{f}
S_{\ge1}(x)-\frac{\tau_{p_1}(x)}{g_{p_1}}
&\le\frac1{g_{p_1}}\sum_{\omega:o\to x}\bigg(\frac{p_1}{|\Lambda|}
 \bigg)^{|\omega|}\sum_{s<t}\underbrace{\prod_{j\ne s,t}\sum_{R_j\in
 \Tcal_{\omega_j}}\bigg(\frac{p_1}{|\Lambda|}\bigg)^{|E_{R_j}|}}_{g_{p_1}
 ^{|\omega|-1}}\nn\\
&\qquad\times\sum_{w\in\Zd}\underbrace{\sum_{R_s\in\Tcal_{\omega_s,w}}
 \bigg(\frac{p_1}{|\Lambda|}\bigg)^{|E_{R_s}|}}_{\tau_{p_1}(w-\omega_s)}\,
 \underbrace{\sum_{R_t\in\Tcal_{\omega_t,w}}\bigg(\frac{p_1}{|\Lambda|}
 \bigg)^{|E_{R_t}|}}_{\tau_{p_1}(\omega_t-w)}\nn\\
&=\frac1{g_{p_1}^2}\sum_{\omega:o\to x}\bigg(\frac1{|\Lambda|}\bigg)^{|\omega|}
 \sum_{s<t}\sum_{w\in\Zd}\tau_{p_1}(w-\omega_s)\,\tau_{p_1}(\omega_t-w)\nn\\
&=\sum_{\substack{y,z\in\Zd\\ (y\ne z)}}\frac{\tau_{p_1}^{*2}(z-y)}{g_{p_1}^2}
 \sum_{\omega:o\to y\to z\to x}\bigg(\frac1{|\Lambda|}\bigg)^{|\omega|}\nn\\
&=\sum_{\substack{y,z\in\Zd\\ (y\ne z)}}\frac{\tau_{p_1}^{*2}(z-y)}{g_{p_1}^2}\,
 S_{\ge0}(y)\,S_{\ge1}(z-y)\,S_{\ge0}(x-z).
\end{align}
The proof of \refeq{RWapprox} is completed by applying 
\refeq{RWapprox-ubd1}--\refeq{RWapprox-ubd2} to 
$(\tau_{p_1}/g_{p_1})^{*2}$ in the above bound.
\QED

\subsection{Detailed analysis of $H$}\label{ss:H}
To complete the proof of Lemma~\ref{lmm:gp1}, it suffices to show the following:

\begin{shaded}
\begin{lmm}\label{lmm:H}
For lattice trees with $d>8$ and $L\uparrow\infty$,
\begin{align}\lbeq{Hlmm}
H=e\sum_{n=3}^\infty\frac{n-1}2D^{*n}(o)+O(\beta^2).
\end{align}
\end{lmm}
\end{shaded}

\Proof{Proof.}
Recall the definition \refeq{Hdef} of $H$:
\begin{align}
H&=\sum_{\substack{Y\subset\Lambda\\ (|Y|\ge2)}}\bigg(\frac{p_1}{|\Lambda|}
 \bigg)^{|Y|}\prod_{y\in Y}\sum_{R_y\in\Tcal_y\setminus\Tcal_o}\bigg(\frac{p_1}
 {|\Lambda|}\bigg)^{|E_{R_y}|}\bigg(1-\prod_{\substack{u,v\in Y\\ (u\ne v)}}
 \ind{V_{R_u}\cap V_{R_v}=\vno}\bigg).
\end{align}
First we split the indicator $1-\prod_{u\ne v}\ind{V_{R_u}\cap V_{R_v}=\vno}$ 
in \refeq{Hdef} by introducing order among pairs of distinct vertices in 
$\Lambda$, called bonds.  If a bond $b$ is earlier than another bond $b'$ in 
that order, we denote it by $b<b'$.  Then we have
\begin{align}\lbeq{Hindicator}
&1-\prod_{\{u,v\}\subset Y}\ind{V_{R_u}\cap V_{R_v}=\vno}\nn\\
&=\sum_{\{u,v\}\subset
 Y}\ind{V_{R_u}\cap V_{R_v}\ne\vno}\prod_{\substack{\{u',v'\}\subset Y\\ (\{u',
 v'\}<\{u,v\})}}\ind{V_{R_{u'}}\cap V_{R_{v'}}=\vno}\nn\\
&=\sum_{\{u,v\}\subset Y}\ind{V_{R_u}\cap V_{R_v}\ne\vno}-\sum_{\{u,v\}\subset
 Y}\ind{V_{R_u}\cap V_{R_v}\ne\vno}\bigg(1-\prod_{\substack{\{u',v'\}\subset Y\\
 (\{u',v'\}<\{u,v\})}}\ind{V_{R_{u'}}\cap V_{R_{v'}}=\vno}\bigg),
\end{align}
where the second sum on the right is zero when $|Y|=2$.  Let $H_1$ be the 
contribution from the first sum on the right:
\begin{align}
H_1=\sum_{\substack{Y\subset\Lambda\\ (|Y|\ge2)}}\bigg(\frac{p_1}{|\Lambda|}
 \bigg)^{|Y|}\prod_{y\in Y}\sum_{R_y\in\Tcal_y\setminus\Tcal_o}\bigg(\frac{p_1}
 {|\Lambda|}\bigg)^{|E_{R_y}|}\sum_{\{u,v\}\subset Y}\ind{V_{R_u}\cap V_{R_v}
 \ne\vno}.
\end{align}
We will later show (after the derivation of \refeq{Hlmm}; see \refeq{Hcomplete}) 
that 
\begin{align}\lbeq{H2bd}
H_2:=H_1-H=O(\beta^2).
\end{align}

Next we investigate $H_1$. Let $H_1'$ be the contribution from the case of $|Y|=2$:
\begin{align}\lbeq{H1'def}
H_1'=\sum_{\{u,v\}\subset\Lambda}\bigg(\frac{p_1}{|\Lambda|}\bigg)^2
 \sum_{\substack{R_u\in\Tcal_u\setminus\Tcal_o \\R_v\in\Tcal_v\setminus\Tcal_o}}
 \bigg(\frac{p_1}{|\Lambda|}\bigg)^{|E_{R_u}|+|E_{R_v}|}\ind{V_{R_u}\cap V_{R_v}
 \ne\vno}.
\end{align}
By subadditivity, we already know that $H_1'=O(\beta)$ for $d>4$.  
By changing the order of sums, we can rewrite 
$H_1-H_1'$ as
\begin{align}
H_1-H_1'&=\sum_{\substack{Y\subset\Lambda\\ (|Y|\ge3)}}\bigg(\frac{p_1}
 {|\Lambda|}\bigg)^{|Y|}\prod_{y\in Y}\sum_{R_y\in\Tcal_y\setminus\Tcal_o}\bigg(
 \frac{p_1}{|\Lambda|}\bigg)^{|E_{R_y}|}\sum_{\{u,v\}\subset Y}\ind{V_{R_u}\cap
 V_{R_v}\ne\vno}\nn \\
&=\underbrace{\sum_{\{u,v\}\subset\Lambda}\bigg(\frac{p_1}{|\Lambda|}\bigg)^2
 \sum_{\substack{R_u\in\Tcal_u\setminus\Tcal_o \\R_v\in\Tcal_v\setminus\Tcal_o}}
 \bigg(\frac{p_1}{|\Lambda|}\bigg)^{|E_{R_u}|+|E_{R_v}|}\ind{V_{R_u}\cap V_{R_v}
 \ne\vno}}_{H_1'}\nn\\
&\qquad\times\sum_{\substack{Y'\subset\Lambda\setminus\{u,v\}\\ (|Y'|\ge1)}}
 \bigg(\frac{p_1}{|\Lambda|}\bigg)^{|Y'|}\prod_{y'\in Y'}\underbrace{\sum_{
 R_{y'}\in\Tcal_{y'}\setminus\Tcal_o}\bigg(\frac{p_1}{|\Lambda|}\bigg)^{|E_{
 R_{y'}}|}}_{g_{p_1}-\tau_{p_1}(y')}.
\end{align}
Similarly to the proof of Lemma~\ref{lmm:G}, we can show that the last line is 
estimated as
\begin{align}
&\sum_{\substack{Y' \subset \Lambda \setminus \{u,v\} \\ (|Y'|\ge 1) }}\bigg(
 \frac1{|\Lambda|}\bigg)^{|Y'|}\prod_{y'\in Y'}\bigg(1-\frac{\tau_{p_1}(y')}
 {g_{p_1}}\bigg)\nn\\
&=\bigg(1+\frac1{|\Lambda|}\bigg)^{|\Lambda|-2}-1+\sum_{\substack{Y'\subset
 \Lambda \setminus \{u,v\} \\ (|Y'|\ge 1) }}\bigg(\frac1{|\Lambda|}\bigg)^{|Y'|}
 \Bigg(\prod_{y'\in Y'}\bigg(1-\frac{\tau_{p_1}(y')}{g_{p_1}}\bigg)-1\Bigg)\nn\\
&=e-1+O(\beta).
\end{align}
Therefore, 
\begin{align}
H_1-H_1'=H_1'\big(e-1+O(\beta)\big),
\end{align}
or equivalently
\begin{align}\lbeq{eH1'}
H_1=eH_1'+O(\beta^2).
\end{align}

Next we investigate $H_1'$.  To do so, we first rewrite 
$\ind{V_{R_u}\cap V_{R_v}\ne\vno}$ in \refeq{H1'def} by introducing order among 
vertices in $\Zd$.  For a vertex set $V$ and an element $x\in V$, we denote by 
$V^{<x}$ the set of vertices in $V$ that are earlier than $x$ in that order.  
By identifying the earliest element $x$ among $V_{R_u}$ that is also in 
$V_{R_v}$ (so that $V^{<x}_{R_u}\cap V_{R_v}=\vno$), we can rewrite 
$\ind{V_{R_u}\cap V_{R_v}\ne\vno}$ as
\begin{align}\lbeq{indintersect}
\ind{V_{R_u}\cap V_{R_v}\ne\vno}&=\sum_{x\in V_{R_u}}\ind{x\in V_{R_v}}\,
 \ind{V^{<x}_{R_u}\cap V_{R_v}=\vno}\nn\\
&=\sum_{x\in\Zd}\ind{x\in V_{R_u}\cap V_{R_v}}-\sum_{x\in V_{R_u}
 \cap V_{R_v}}\Big(1-\ind{V^{<x}_{R_{u}} \cap V_{R_{v}}= \vno}\Big).
\end{align}
Let $H_1''$ be the contribution from the first sum in the last line:
\begin{align}\lbeq{H1''def}
H_1''&=\sum_{\{u,v\}\subset\Lambda}\bigg(\frac{p_1}{|\Lambda|}\bigg)^2\sum_{x
 \in\Zd}\sum_{\substack{R_u\in\Tcal_{u,x}\setminus\Tcal_o \\R_v\in\Tcal_{v,x}
 \setminus\Tcal_o}}\bigg(\frac{p_1}{|\Lambda|}\bigg)^{|E_{R_u}|+|E_{R_v}|}\nn\\
&=\sum_{\{u,v\}\subset\Lambda}\bigg(\frac{p_1}{|\Lambda|}\bigg)^2\sum_{x\in\Zd}
 \Big(\tau_{p_1}(u-x)-\tau_{p_1}^{\sss(3)}(o,u,x)\Big)\Big(\tau_{p_1}(v-x)
 -\tau_{p_1}^{\sss(3)}(o,v,x)\Big)\nn\\
&\stackrel{\text{\refeq{p1gp1}}}=\sum_{\{u,v\}\subset\Lambda}\bigg(\frac1
 {|\Lambda|}\bigg)^2\sum_{x\in\Zd}\bigg(\frac{\tau_{p_1}(u-x)}{g_{p_1}}
 -\frac{\tau_{p_1}^{\sss(3)}(o,u,x)}{g_{p_1}}\bigg)\bigg(\frac{\tau_{p_1}(v-x)}
 {g_{p_1}}-\frac{\tau_{p_1}^{\sss(3)}(o,v,x)}{g_{p_1}}\bigg),
\end{align}
where $\tau_p^{\sss(3)}(o,u,x)$ is a 3-point function, defined as
\begin{align}
\tau_p^{\sss(3)}(o,u,x)=\sum_{T\in \Tcal_{o,u,x}}W_p(T).
\end{align}
We will later show that 
\begin{align}\lbeq{H2''bd}
H_2''&:=H_1''-H_1'=O(\beta^2).
\end{align}

Finally we investigate $H_1''$.  The dominant contribution to $H_1''$ comes 
from the product of 2-point functions:
\begin{align}\lbeq{H1'''def}
H_1'''&:=\sum_{\{u,v\}\subset\Lambda}\bigg(\frac1{|\Lambda|}\bigg)^2\sum_{x\in\Zd}
 \frac{\tau_{p_1}(u-x)}{g_{p_1}}\frac{\tau_{p_1}(v-x)}{g_{p_1}}\nn\\
&=2\sum_{\{u,v\}\subset\Lambda}\bigg(\frac1{|\Lambda|}\bigg)^2\frac{\tau_{p_1}
 (u-v)}{g_{p_1}}+\sum_{\{u,v\}\subset\Lambda}\bigg(\frac1{|\Lambda|}\bigg)^2
 \sum_{x\ne u,v}\frac{\tau_{p_1}(u-x)}{g_{p_1}}\frac{\tau_{p_1}(v-x)}{g_{p_1}},
\end{align}
where we have used the identity $\tau_p(o)=g_p$.  We will later show that the other
contribution to $H_1''$ which involves 3-point functions is estimated as 
\begin{align}\lbeq{H2'''bd}
H_2''':=H_1''-H_1'''=O(\beta^2).
\end{align}
By Lemma~\ref{lmm:RWapprox}, the first term in \refeq{H1'''def}
is estimated as 
\begin{align}\lbeq{H1''main1}
&2\sum_{\{u,v\}\subset\Lambda}\bigg(\frac1{|\Lambda|}\bigg)^2\frac{\tau_{p_1}
 (u-v)}{g_{p_1}}\nn\\
&=\sum_{\substack{u,v\in\Lambda\\ (u\ne v)}}\bigg(\frac1{|\Lambda|}\bigg)^2
 S_{\ge1}(u-v)+\sum_{\substack{u,v\in\Lambda\\ (u\ne v)}}\bigg(\frac1{|\Lambda|}\bigg)^2\left(\frac{\tau_{p_1}
 (u-v)}{g_{p_1}}- S_{\ge1}(u-v)\right)
 \nn\\
&=S_{\ge3}(o)-\underbrace{\frac1{|\Lambda|}S_{\ge1}(o)}_{O(\beta^2)\text{ for }
 d>2}+\underbrace{\sum_{\substack{u,v\in\Lambda\\ (u\ne v)}}\bigg(\frac1{|\Lambda|}\bigg)^2\left(\frac{\tau_{p_1}
 (u-v)}{g_{p_1}}- S_{\ge1}(u-v)\right)}_{\overset{\refeq{RWapprox}-\refeq{RWapproxdiagramA}}{\le} \raisebox{-1pc}{\includegraphics[scale=0.3]{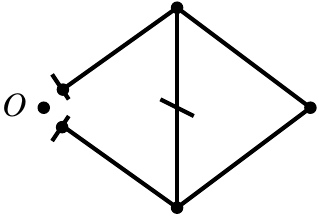}}}
\end{align}
where a gap next to the origin in the last diagram represents $1/|\Lambda|$. By translation invariance and \refeq{RWanal1}-\refeq{RWanal2}, the last term is bounded above by
\begin{align}
\raisebox{-1pc}{\includegraphics[scale=0.3]{H1ppp1e}}=\raisebox{-1pc}{\includegraphics[scale=0.3]{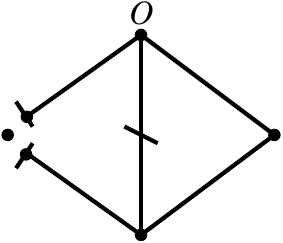}}=\sum_{\substack{y \in \Zd \\ y\neq o}}\underbrace{\raisebox{-1.5pc}{\includegraphics[scale=0.3]{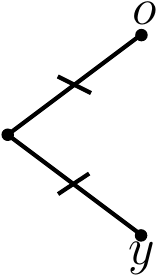}}}_{\le \|S_{\ge1}^{*2}\|_\infty}~~\raisebox{-1.5pc}{\includegraphics[scale=0.3]{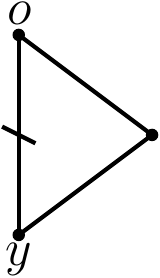}}
&\le \|S_{\ge1}^{*2}\|_\infty(S_{\ge0}^{*2}*S_{\ge1})(o)=O(\beta^2).
\end{align}
Similarly, the second term in \refeq{H1'''def} is estimated as
\begin{align}\lbeq{H1''main2}
&\sum_{\{u,v\}\subset\Lambda}\bigg(\frac1{|\Lambda|}\bigg)^2\sum_{x\ne u,v}
 \frac{\tau_{p_1}(u-x)}{g_{p_1}}\frac{\tau_{p_1}(v-x)}{g_{p_1}}\nn\\
&=\frac12S_{\ge2}^{*2}(o)-\underbrace{\frac1{2|\Lambda|}S_{\ge1}^{*2}(o)}_{O
 (\beta^2)\text{ for }d>4}\nn\\
&+\frac12\sum_{\substack{u,v\in\Lambda \\ u\neq v}}\bigg(\frac1{|\Lambda|}\bigg)^2\sum_{x\ne u,v}\underbrace{\left(\frac{\tau_{p_1}(u-x)}{g_{p_1}}\frac{\tau_{p_1}(v-x)}{g_{p_1}}-S_{\ge1}(u-x)S_{\ge1}(v-x)\right)}_{\le~ 2~\raisebox{-1pc}{\includegraphics[scale=0.3]{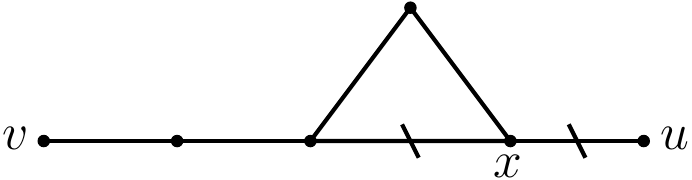}}~+~ \raisebox{-1pc}{\includegraphics[scale=0.3]{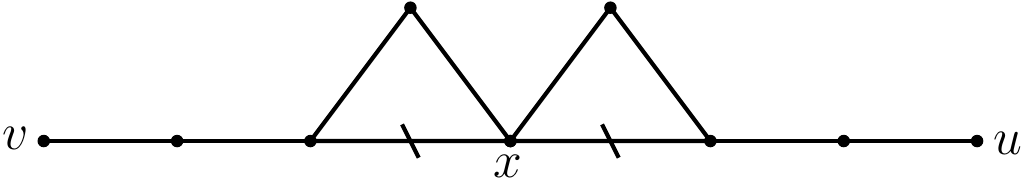}}}.
\end{align}
By Lemma~\ref{lmm:RWapprox}, \refeq{RWapproxdiagramB} and the translation invariance, the last term of \refeq{H1''main2} is bounded above by
\begin{align}\lbeq{d}
\raisebox{-1pc}{\includegraphics[scale=0.3]{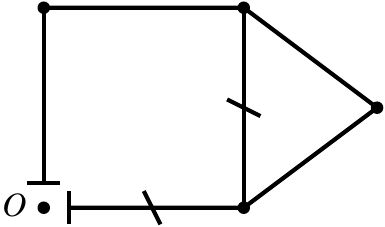}}
 +\frac12\raisebox{-1.2pc}{\includegraphics[scale=0.28]{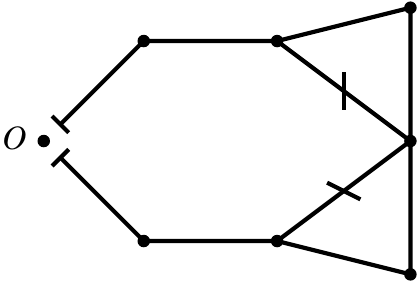}}&=\sum_{\substack{y \in \Zd \\ y \neq o}}\underbrace{\raisebox{-1.5pc}{\includegraphics[scale=0.3]{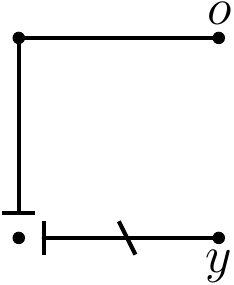}}}_{\le \|S_{\ge1}^{*3}\|_\infty}~~\raisebox{-1.4pc}{\includegraphics[scale=0.3]{H1ppp1eE}}+\frac12\sum_{\substack{y,z \in \Zd \\ y,z \neq o}}\underbrace{\raisebox{-1.5pc}{\includegraphics[scale=0.3]{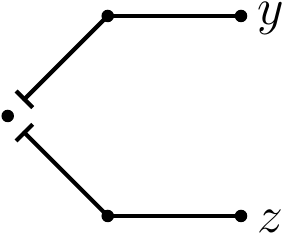}}}_{\le \|S_{\ge0}^{*2}*S_{\ge1}^{*2}\|_\infty}~~\raisebox{-1.8pc}{\includegraphics[scale=0.3]{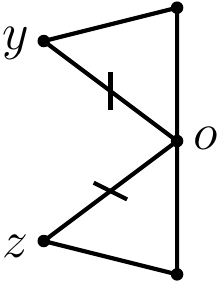}}\nn\\
&\le \underbrace{\|S_{\ge1}^{*3}\|_\infty(S_{\ge0}^{*2}*S_{\ge1})(o)}_{O(\beta^2)
 \text{ for }d>6}+\frac12\underbrace{\|S_{\ge0}^{*2}*S_{\ge1}^{*2}\|_\infty
 (S_{\ge0}^{*2}*S_{\ge1})(o)^2}_{O(\beta^3)\text{ for }d>8}\nn\\
&=O(\beta^2).
\end{align}
Therefore, 
\begin{align}\lbeq{H_1'''bd}
H_1'''=S_{\ge3}(o)+\frac12S_{\ge2}^{*2}(o)+O(\beta^2)&=\sum_{n=3}^\infty D^{*n}
 (o)+\frac12\sum_{n,m=2}^\infty D^{*(n+m)}(o)+O(\beta^2)\nn\\
&=\sum_{n=3}^\infty\frac{n-1}2D^{*n}(o)+O(\beta^2).
\end{align}

Summarizing all the above estimates, we arrive at
\begin{align}\lbeq{Hcomplete}
H\stackrel{\text{\refeq{H2bd}}}=H_1+O(\beta^2)
 \stackrel{\text{\refeq{eH1'}}}=eH_1'+O(\beta^2)
 \stackrel{\text{\refeq{H2''bd}}}=eH_1''+O(\beta^2)
 \stackrel{\text{\refeq{H2'''bd}}}=eH_1'''+O(\beta^2)\nn\\
\stackrel{\text{\refeq{H_1'''bd}}}=e\sum_{n=3}^\infty\frac{n-1}2D^{*n}(o)
 +O(\beta^2),
\end{align}
as required.  It remains to show \refeq{H2bd}, \refeq{H2''bd} and 
\refeq{H2'''bd}.
\QED

\Proof{Proof of \refeq{H2bd}: bounding $H_2$.}
First we recall that $H_2$ is the contribution from the second sum on the right 
of \refeq{Hindicator}:
\begin{align}\lbeq{H2def}
H_2=\sum_{\substack{Y\subset\Lambda\\ (|Y|\ge3)}}\bigg(\frac{p_1}{|\Lambda|}
 \bigg)^{|Y|}\prod_{y\in Y}\sum_{R_y\in\Tcal_y\setminus\Tcal_o}\bigg(\frac{p_1}
 {|\Lambda|}\bigg)^{|E_{R_y}|}\sum_{\{u,v\}\subset Y}\ind{V_{R_u}\cap V_{R_v}
 \ne\vno}\nn\\
\times\bigg(1-\prod_{\substack{\{u',v'\}\subset Y\\ (\{u',v'\}<\{u,v\})}}
 \ind{V_{R_{u'}}\cap V_{R_{v'}}=\vno}\bigg),
\end{align}
which is nonnegative. Since we get an upper bound 
\begin{align}\lbeq{e}
1-\prod_{\substack{\{u',v'\}\subset Y\\ (\{u',v'\}<\{u,v\})}}
 \ind{V_{R_{u'}}\cap V_{R_{v'}}=\vno}\le \sum_{\substack{ \{u',v'\}
 \subset Y\\ (\{u',v'\}<\{u,v\})}}\ind{V_{R_{u'}}\cap V_{R_{v'}}\ne\vno}
\end{align}
in a same manner as \refeq{Hindicator},
we can bound $H_2$ as
\begin{align}\lbeq{H2bd-pr1}
H_2&\le\sum_{\substack{Y\subset\Lambda\\ (|Y|\ge3)}}\bigg(\frac{p_1}{|\Lambda|}
 \bigg)^{|Y|}\prod_{y\in Y}\sum_{R_y\in\Tcal_y\setminus\Tcal_o}\bigg(\frac{p_1}
 {|\Lambda|}\bigg)^{|E_{R_y}|}\sum_{\substack{\{u,v\}\subset Y\\ \{u',v'\}
 \subset Y\\ (\{u',v'\}<\{u,v\})}}\ind{V_{R_u}\cap V_{R_v}\ne\vno}\,
 \ind{V_{R_{u'}}\cap V_{R_{v'}}\ne\vno}\nn\\
&=\frac{1}{2}\sum_{\substack{Y\subset\Lambda\\ (|Y|\ge3)}}\bigg(\frac{p_1}{|\Lambda|}
 \bigg)^{|Y|}\prod_{y\in Y}\sum_{R_y\in\Tcal_y\setminus\Tcal_o}\bigg(\frac{p_1}
 {|\Lambda|}\bigg)^{|E_{R_y}|}\sum_{\substack{\{u,v\}\subset Y\\ \{u',v'\}
 \subset Y\\ (\{u',v'\}\neq\{u,v\})}}\ind{V_{R_u}\cap V_{R_v}\ne\vno}\,
 \ind{V_{R_{u'}}\cap V_{R_{v'}}\ne\vno}.
\end{align}
Since $\{u,v\}\ne\{u',v'\}$, the union $\{u,v\}\cup\{u',v'\}$ consists of 
either three or four distinct vertices.  We denote the contribution from the 
former by $H_{2,3}$, and that from the latter by $H_{2,4}$ and then we obtain
\begin{align}\lbeq{H2<1/2H23+H24}
H_2\le \frac{1}{2}(H_{2,3}+H_{2,4}).
\end{align}

\begin{figure}[t]
\begin{center}
\includegraphics[scale=0.6]{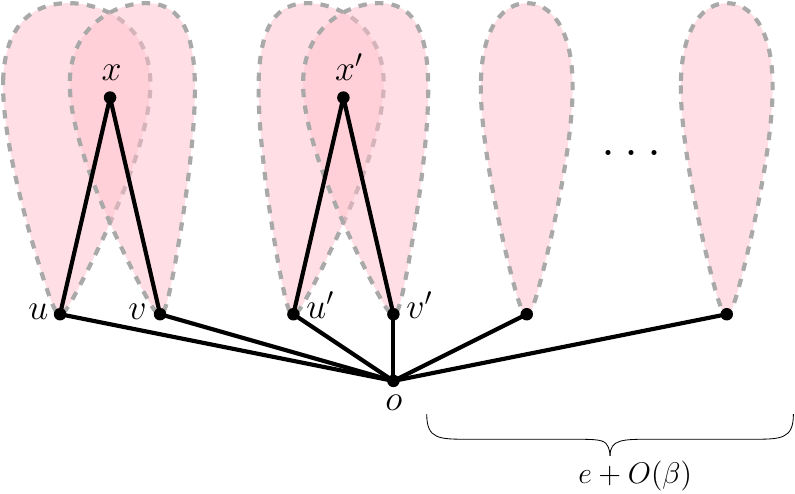}
\end{center}
\caption{Schematic representation of $H_{2,4}$.}
\label{fig:H24}
\end{figure}

First we investigate $H_{2,4}$, which is bounded as (see Figure~\ref{fig:H24})
\begin{align}\lbeq{H2bd-pr2}
H_{2,4}&=\sum_{\substack{Y\subset\Lambda\\ (|Y|\ge4)}}\bigg(\frac{p_1}{|\Lambda|}
 \bigg)^{|Y|}\prod_{y\in Y}\sum_{R_y\in\Tcal_y\setminus\Tcal_o}\bigg(\frac{p_1}
 {|\Lambda|}\bigg)^{|E_{R_y}|}\sum_{\substack{u,v,u',v'\in Y\\
 (\text{distinct})}}\ind{V_{R_u}\cap V_{R_v}\ne\vno}\ind{V_{R_{u'}}\cap
 V_{R_{v'}}\ne\vno}\nn\\
&=\sum_{\substack{u,v,u',v'\in\Lambda\\ (\text{distinct})}}\bigg(\frac{p_1}
 {|\Lambda|}\bigg)^4\!\!\sum_{\substack{R_u\in\Tcal_u\setminus\Tcal_o\\ R_v\in
 \Tcal_v\setminus\Tcal_o\\ R_{u'}\in\Tcal_{u'}\setminus\Tcal_o\\ R_{v'}\in
 \Tcal_{v'}\setminus\Tcal_o}}\!\!\bigg(\frac{p_1}{|\Lambda|}\bigg)^{|E_{R_u}|
 +|E_{R_v}|+|E_{R_{u'}}|+|E_{R_{v'}}|}\ind{V_{R_u}\cap V_{R_v}\ne\vno}
 \ind{V_{R_{u'}}\cap V_{R_{v'}}\ne\vno}\nn\\
&\qquad\times\Bigg(1+\sum_{\substack{Y'\subset\Lambda\setminus\{u,v,u',v'\}\\
 (|Y'|\ge1)}}\bigg(\frac{p_1}{|\Lambda|}\bigg)^{|Y'|}\prod_{y\in Y'}\underbrace{
 \sum_{R_y\in\Tcal_y\setminus\Tcal_o}\bigg(\frac{p_1}{|\Lambda|}\bigg)^{|E_{
 R_y}|}}_{g_{p_1}-\tau_{p_1}(y)}\Bigg).
\end{align}
By \refeq{p1gp1} (i.e., $p_1g_{p_1}=1$), the last line is equal to
\begin{align}
1+\sum_{\substack{Y'\subset\Lambda\setminus\{u,v,u',v'\}\\ (|Y'|\ge1)}}\bigg(
 \frac1{|\Lambda|}\bigg)^{|Y'|}=\bigg(1+\frac1{|\Lambda|}\bigg)^{|\Lambda|-4}
 =e+O(\beta).
\end{align}
Then, by ignoring the constraint that $\{u,v\}$ and $\{u',v'\}$ are disjoint 
pairs and using the trivial inequality 
$\ind{V_{R_u}\cap V_{R_v}\ne\vno}\le\sum_x\ind{x\in V_{R_u}\cap V_{R_v}}$ 
as well as the relation $\Tcal_{u,x}\setminus\Tcal_o\subset\Tcal_{u,x}$,
$H_{2,4}$ is further bounded above as
\begin{align}\lbeq{H2bd-pr3}
H_{2,4}&\le\big(e+O(\beta)\big)\Bigg(\sum_{\{u,v\}\subset\Lambda}
 \bigg(\frac{p_1}{|\Lambda|}\bigg)^2\sum_{\substack{R_u\in\Tcal_u\setminus
 \Tcal_o\\ R_v\in\Tcal_v\setminus\Tcal_o}}\bigg(\frac{p_1}{|\Lambda|}
 \bigg)^{|E_{R_u}|+|E_{R_v}|}\ind{V_{R_u}\cap V_{R_v}\ne\vno}\Bigg)^2\nn\\
&\le\big(e+O(\beta)\big)\Bigg(\sum_{\{u,v\}\subset\Lambda}\bigg(\frac{p_1}
 {|\Lambda|}\bigg)^2\sum_{x\in\Zd}\underbrace{\sum_{\substack{R_u\in
 \Tcal_{u,x}\\ R_v\in\Tcal_{v,x}}}\bigg(\frac{p_1}{|\Lambda|}\bigg)^{|E_{R_u}|
 +|E_{R_v}|}}_{\tau_{p_1}(u-x)\,\tau_{p_1}(v-x)}\Bigg)^2\nn\\
&\stackrel{\text{\refeq{p1gp1}}}=\big(e+O(\beta)\big)\Bigg(\sum_{\{u,v\}\subset
 \Lambda}\bigg(\frac1{|\Lambda|}\bigg)^2\sum_{x\in\Zd}\frac{\tau_{p_1}(u-x)}
 {g_{p_1}}\frac{\tau_{p_1}(v-x)}{g_{p_1}}\Bigg)^2.
\end{align}
Finally, by using $\tau_{p_1}(u-x)/g_{p_1}\le S_{\ge1}(u-x)$ for $x\ne u$ 
(cf., the first inequality in \refeq{RWapprox}) and $\tau_p(o)=g_p$ for $x=u$, 
we arrive at
\begin{align}\lbeq{H2bd-pr31}
H_{2,4}&\le\big(e+O(\beta)\big)\Bigg(\sum_{u,v\in\Lambda}\bigg(\frac1
 {|\Lambda|}\bigg)^2S_{\ge0}^{*2}(u-v)\Bigg)^2\le\big(e+O(\beta)\big)\,S_{\ge
 1}^{*2}(o)^2\,\stackrel{d>4}=O(\beta^2).
\end{align}

\begin{figure}[t]
\begin{center}
\includegraphics[scale=0.6]{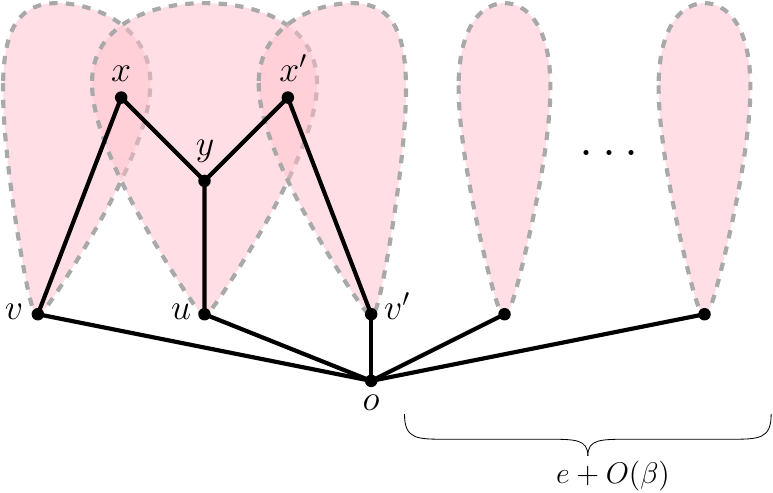}
\end{center}
\caption{Schematic representation of $H_{2,3}$.}
\label{fig:H23}
\end{figure}

Next we investigate $H_{2,3}$, which is bounded in a similar way to 
\refeq{H2bd-pr3} as (see Figure~\ref{fig:H23})
\begin{align}\lbeq{H2bd-pr4}
H_{2,3}&=\sum_{\substack{Y\subset\Lambda\\ (|Y|\ge3)}}\bigg(\frac{p_1}{|\Lambda|}
 \bigg)^{|Y|}\prod_{y\in Y}\sum_{R_y\in\Tcal_y\setminus\Tcal_o}\bigg(\frac{p_1}
 {|\Lambda|}\bigg)^{|E_{R_y}|}\sum_{\substack{u,v,v'\in Y\\ (\text{distinct})}}
 \ind{V_{R_u}\cap V_{R_v}\ne\vno}\,\ind{V_{R_u}\cap V_{R_{v'}}\ne\vno}\nn\\
&\le\big(e+O(\beta)\big)\sum_{\substack{u,v,v'\in\Lambda\\ (\text{distinct})}}
 \bigg(\frac{p_1}{|\Lambda|}\bigg)^3\sum_{x,x'\in\Zd}\sum_{
 \substack{R_u\in\Tcal_{u,x,x'}\\ R_v\in\Tcal_{v,x}\\ R_{v'}\in\Tcal_{v',x'}}}
 \bigg(\frac{p_1}{|\Lambda|}\bigg)^{|E_{R_u}|+|E_{R_v}|+|E_{R_{v'}}|}\nn\\
&=\big(e+O(\beta)\big)g_{p_1}^2\sum_{\substack{u,v,v'\in\Lambda\\
 (\text{distinct})}}\bigg(\frac1{|\Lambda|}\bigg)^3\sum_{x,x'\in\Zd}
 \frac{\tau_{p_1}^{\sss(3)}(u,x,x')}{g_{p_1}^3}\frac{\tau_{p_1}(v-x)}{g_{p_1}}
 \frac{\tau_{p_1}(v'-x')}{g_{p_1}}.
\end{align}
Due to submultiplicativity, we can bound $\tau_p^{\sss(3)}(u,x,x')$ as 
\begin{align}\lbeq{bdof3pf}
\tau_p^{\sss(3)}(u,x,x')\le\sum_{y\in\Zd}\tau_p(u-y)\,\tau_p(x-y)\,\tau_p(x'-y).
\end{align}
Then, by using $\tau_{p_1}(u-x)/g_{p_1}\le S_{\ge1}(u-x)$ for $x\ne u$ and 
$\tau_p(o)=g_p$ for $x=u$, we can bound the sum in \refeq{H2bd-pr4} as
\begin{align}\lbeq{H2bd-pr5}
&\sum_{\substack{u,v,v'\in\Lambda\\ (\text{distinct})}}\bigg(\frac1
 {|\Lambda|}\bigg)^3\sum_{y,x,x'\in\Zd}\frac{\tau_{p_1}(u-y)}{g_{p_1}}
 \frac{\tau_{p_1}(x-y)}{g_{p_1}}\frac{\tau_{p_1}(x'-y)}{g_{p_1}}\frac{\tau_{p_1}
 (v-x)}{g_{p_1}}\frac{\tau_{p_1}(v'-x')}{g_{p_1}}\nn\\
&\le\sum_{y,x,x'\in\Zd}S_{\ge0}(x-y)\,S_{\ge0}(x'-y)\underbrace{
 \sum_{u,v,v'\in\Lambda}\bigg(\frac1{|\Lambda|}\bigg)^3S_{\ge0}(u-y)\,S_{\ge
 0}(v-x)\,S_{\ge0}(v'-x')}_{S_{\ge1}(y)\,S_{\ge1}(x)\,S_{\ge1}(x')}\nn\\
&=\sum_{y\in\Zd}(S_{\ge0}*S_{\ge1})(y)^2~S_{\ge1}(y)\nn\\
&\le\|S_{\ge0}*S_{\ge1}\|_\infty\,(S_{\ge0}*S_{\ge1}^{*2})(o)\stackrel{d>6}=
 O(\beta^2).
\end{align}
This together with \refeq{H2<1/2H23+H24} and \refeq{H2bd-pr31} implies
\begin{align}\lbeq{H2<H23+H24}
H_2\le \frac{1}{2}(H_{2,3}+H_{2,4})=O(\beta^2),
\end{align}
as required.
\QED

\begin{figure}[t]
\begin{center}
\includegraphics[scale=0.6]{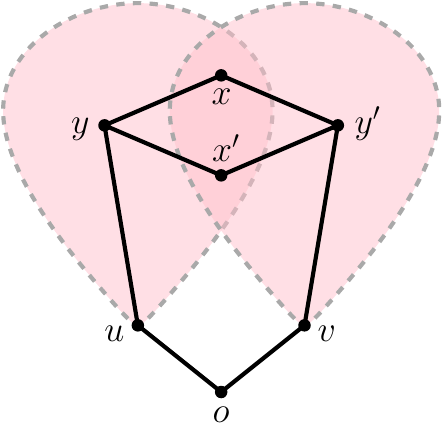}
\end{center}
\caption{Schematic representation of the bound on $H_2''$ due to \refeq{H7ind}.}
\label{fig:H2ppbd}
\end{figure}

\Proof{Proof of \refeq{H2''bd}: bounding $H_2''$.}
First we recall that $H_2''$ is the contribution to $H_1'$ from the second sum 
on the right of \refeq{indintersect}:
\begin{align}
H_2''=\sum_{\{u,v\}\subset\Lambda}\bigg(\frac{p_1}{|\Lambda|}\bigg)^2\sum_{x
 \in\Zd}\sum_{\substack{R_u\in\Tcal_{u,x}\setminus\Tcal_o \\R_v\in\Tcal_{v,x}
 \setminus\Tcal_o}}\bigg(\frac{p_1}{|\Lambda|}\bigg)^{|E_{R_u}|+|E_{R_v}|}\Big(1
 -\ind{V^{<x}_{R_{u}} \cap V_{R_{v}}= \vno}\Big).
\end{align}
Notice that
\begin{align}\lbeq{H7ind}
1-\ind{V^{<x}_{R_{u}} \cap V_{R_v}=\vno}\le\sum_{x'\in\Zd\setminus\{x\}}
 \ind{x'\in V_{R_u}\cap V_{R_v}}.
\end{align}
By the inclusion relation $\Tcal_{u,x,x'}\setminus\Tcal_o\subset\Tcal_{u,x,x'}$ 
and using \refeq{bdof3pf}, \refeq{p1gp1} and \refeq{RWapprox}, we can bound 
$H_2''$ as (see Figure~\ref{fig:H2ppbd})
\begin{align}\lbeq{includexw}
H_2''&\le\sum_{\{u,v\}\subset\Lambda}\bigg(\frac{p_1}{|\Lambda|}\bigg)^2
 \sum_{\substack{x,x'\in\Zd\\ (x\ne x')}}\underbrace{\sum_{\substack{R_u\in
 \Tcal_{u,x,x'}\\ R_v\in\Tcal_{v,x,x'}}}\bigg(\frac{p_1}{|\Lambda|}\bigg)^{
 |E_{R_u}|+|E_{R_v}|}}_{\tau_{p_1}^{\sss(3)}(u,x,x')\,\tau_{p_1}^{\sss(3)}(v,x,
 x')}\nn\\
&\le g_{p_1}^4\sum_{\substack{x,x',y,y'\in\Zd\\ (x\ne x')}}\frac{\tau_{p_1}(y
 -x)}{g_{p_1}}\frac{\tau_{p_1}(y-x')}{g_{p_1}}\frac{\tau_{p_1}(y'-x)}{g_{p_1}}
 \frac{\tau_{p_1}(y'-x')}{g_{p_1}}\nn\\
&\qquad\times\frac12\sum_{u,v\in\Lambda}\bigg(\frac1{|\Lambda|}\bigg)^2
 \frac{\tau_{p_1}(u-y)}{g_{p_1}}\frac{\tau_{p_1}(v-y')}{g_{p_1}}\nn\\
&\le\frac{g_{p_1}^4}2\!\!\sum_{\substack{x,x',y,y'\in\Zd\\ (x\ne x')}}\!\!S_{\ge
 0}(y-x)\,S_{\ge0}(y-x')\,S_{\ge0}(y'-x)\,S_{\ge0}(y'-x')\,S_{\ge1}(y)\,S_{\ge1}
 (y').
\end{align}
Shifting the variables by $-x'$ and changing the variables $x-x',y-x',y'-x'$ to 
the new ones $w,z,z'$, respectively, we can rewrite the above sum as
\begin{align}
\sum_{\substack{w,z,z'\in\Zd\\ (w\ne o)}}S_{\ge0}(z-w)\,S_{\ge0}(z)\,S_{\ge0}
 (z'-w)\,S_{\ge0}(z')\underbrace{\sum_{x'\in\Zd}S_{\ge1}(z+x')\,S_{\ge1}(z'
 +x')}_{S_{\ge1}^{*2}(z-z')},
\end{align}
which is bounded above by
\begin{align}
\|S_{\ge1}^{*2}\|_\infty\sum_{w\ne o}S_{\ge0}^{*2}(w)^2
&=\|S_{\ge1}^{*2}\|_\infty\sum_{w\ne o}\bigg(\sum_{n=1}^\infty(n+1)D^{*n}(w)
 \bigg)^2\nn\\
&\le\|S_{\ge1}^{*2}\|_\infty\sum_{t=2}^\infty D^{*t}(o)\underbrace{\sum_{n=1}^{t
 -1}(n+1)(t-n+1)}_{O(t^3)}~\stackrel{d>8}=O(\beta^2),
\end{align}
as required.
\QED

\Proof{Proof of \refeq{H2'''bd}: bounding $H_2'''$.}
First we recall that $H_2'''$ is the contribution to $H_1''$ which involves 
3-point functions (cf., \refeq{H1''def}):
\begin{align}
H_2'''&=\sum_{\{u,v\}\subset\Lambda}\bigg(\frac1{|\Lambda|}\bigg)^2\sum_{x\in
 \Zd}\frac{\tau_{p_1}^{\sss(3)}(o,u,x)}{g_{p_1}}\bigg(\frac{\tau_{p_1}^{\sss(3)}
 (o,v,x)}{g_{p_1}}-2\frac{\tau_{p_1}(v-x)}{g_{p_1}}\bigg).
\end{align}
By \refeq{p1gp1}, \refeq{bdof3pf} and \refeq{RWapprox}, we can readily conclude 
that
\begin{align}
|H_2'''|&\le\frac{g_{p_1}^2}2\sum_{u,v\in\Lambda}\bigg(\frac1{|\Lambda|}\bigg)^2
 \sum_{x,y\in\Zd}\frac{\tau_{p_1}(y)}{g_{p_1}}\frac{\tau_{p_1}(y-u)}
 {g_{p_1}}\frac{\tau_{p_1}(y-x)}{g_{p_1}}\nn\\
&\qquad\times\bigg(g_{p_1}^2\sum_{z\in\Zd}\frac{\tau_{p_1}(z)}{g_{p_1}}
 \frac{\tau_{p_1}(z-v)}{g_{p_1}}\frac{\tau_{p_1}(z-x)}{g_{p_1}}
 +2\frac{\tau_{p_1}(v-x)}{g_{p_1}}\bigg)\nn\\
&\le g_{p_1}^2\sum_{x,y\in\Zd}S_{\ge0}(y)\,S_{\ge1}(y)\,S_{\ge0}(y-x)\nn\\
&\qquad\times\bigg(\frac{g_{p_1}^2}2\sum_{z\in\Zd}S_{\ge0}(z)\,S_{\ge1}(z)\,
 S_{\ge0}(z-x)+S_{\ge1}(x)\bigg)\nn\\
&=g_{p_1}^2\sum_{y\in\Zd}S_{\ge0}(y)\,S_{\ge1}(y)\bigg(\frac{g_{p_1}^2}2
 \sum_{z\in\Zd}S_{\ge0}(z)\,S_{\ge1}(z)\,S_{\ge0}^{*2}(y-z)+(S_{\ge0}*S_{\ge1})
 (y)\bigg)\nn\\
&\le g_{p_1}^2\underbrace{(S_{\ge0}*S_{\ge1})(o)}_{O(\beta)\text{ for }d>4}
 \bigg(\frac{g_{p_1}^2}2(S_{\ge0}*S_{\ge1})(o)\underbrace{\|S_{\ge0}^{*2}
 \|_\infty}_{O(1)\text{ for }d>4}+\underbrace{\|S_{\ge0}*S_{\ge1}\|_\infty}_{O
 (\beta)\text{ for }d>4}\bigg)\nn\\
&=O(\beta^2),
\end{align}
as required.
\QED

\section{Difference between lattice trees and lattice animals}\label{s:difference}
Finally we prove Theorem~\ref{thm:main} for lattice animals.  
Recall that, by Lemma~\ref{lmm:pc-2nd}, it suffices to investigate 
$p_1=1/g_{p_1}$ (cf., \refeq{p1gp1}).  The following is the key lemma:

\begin{shaded}
\begin{lmm}\label{lmm:gp1a}
For lattice animals with $d>8$ and $L\uparrow\infty$,
\begin{align}\lbeq{animmalgp1}
g_{p_1}=e\bigg(1-\sum_{n=2}^\infty\frac{n+1}2D^{*n}(o)\bigg)
 +\frac12S_{\ge3}(o)+O(\beta^2).
\end{align}
Consequently,
\begin{align}
p_1=\frac1e+\sum_{n=2}^\infty\frac{n+1}{2e}D^{*n}(o)-\frac1{2e^2}S_{\ge3}(o)+O(\beta^2).
\end{align}
\end{lmm}
\end{shaded}

\Proof{Proof.}
As a first step, we want a similar decomposition to \refeq{gp1-rewr} for 
lattice animals.  To do so, we identify the connected neighbors $Y$ of the 
origin, just as done in \refeq{gp1-rewr}.  Then, we introduce $\Gamma(Y)$, 
which is the set of all partitions of $Y$.  For example, if $Y=\{1,2,3\}$, then 
\begin{align}
\Gamma(Y)=\bigg\{\{Y\},~\big\{\{1,2\},\{3\}\big\},~\big\{\{1,3\},\{2\}\big\},~
 \big\{\{1\},\{2,3\}\big\},~\big\{\{1\},\{2\},\{3\}\big\}\bigg\}.
\end{align}
For a partition $\gamma\in\Gamma(Y)$, we denote by $|\gamma|$ the number of 
sets in $\gamma$, so that $\gamma=\{\gamma_j\}_{j=1}^{|\gamma|}$.  We can 
rewrite $g_{p_1}$ as 
\begin{align}
g_{p_1}&=1+\sum_{\substack{Y\subset\Lambda\\ (|Y|\ge1)}}\bigg(\frac{p_1}
 {|\Lambda|}\bigg)^{|Y|}\sum_{\gamma\in\Gamma(Y)}\,\prod_{j=1}^{|\gamma|}\,
 \sum_{R_j\in\Acal_{\gamma_j}\setminus\Acal_o}\bigg(\frac{p_1}{|\Lambda|}
 \bigg)^{|E_{R_j}|}\prod_{i<j}\ind{V_{R_i}\cap V_{R_j}=\vno}.
\end{align} 
The contribution from the maximum partition 
$\bar\gamma=\{\{y\}\}_{y\in Y}$ (i.e., $|\bar\gamma|=|Y|$) is equal to
\refeq{gp1-rewr} (with $\Tcal$ replaced by $\Acal$) and can 
be decomposed into $G$ and $H$ as in \refeq{gp1=G-H} (with $R_y$ regarded as 
animals instead of trees).  Let $I$ be the contribution from the remaining 
partitions $\gamma\in\Gamma(Y)$ with $|\gamma|<|Y|$, which is zero for 
lattice trees:
\begin{align}\lbeq{LTLAdiff}
I&=g_{p_1}-(G-H)\nn\\
&=\sum_{\substack{Y\subset\Lambda\\ (|Y|\ge2)}}\bigg(\frac{p_1}{|\Lambda|}
 \bigg)^{|Y|}\sum_{\substack{\gamma\in\Gamma(Y)\\ (|\gamma|<|Y|)}}\prod_{j
 =1}^{|\gamma|}\sum_{R_j\in\Acal_{\gamma_j}\setminus\Acal_o}\bigg(\frac{p_1}
 {|\Lambda|}\bigg)^{|E_{R_j}|}\prod_{i<j}\ind{V_{R_i}\cap V_{R_j}=\vno}.
\end{align}

To evaluate $G, H$ and $I$ for lattice animals, we cannot apply 
Lemma~\ref{lmm:RWapprox}, which is a powerful tool for lattice trees to 
identify the coefficients of $\beta$ as well as to estimate the error terms of 
$O(\beta^2)$.  For the latter purpose for lattice animals, we will use the 
infrared bound \refeq{MFbehavior} (and monotonicity in $p$, i.e., 
$\tau_{p_1}\le\tau_{\pc}$); for the former purpose, we will use the 
following bounds that correspond to \refeq{RWest1LT}, \refeq{H1''main1} and 
\refeq{H1''main2}, respectively:

\begin{shaded}
\begin{lmm}\label{lmm:RWapproxLA}
For lattice animals with $d>8$ and $L\uparrow\infty$,
\begin{gather}
\bigg|\sum_{u\in\Lambda}\frac1{|\Lambda|}\frac{\tau_{p_1}(u)}{g_{p_1}}
 -S_{\ge2}(o)\bigg|=O(\beta^2),\lbeq{RWest1LA}\\
\bigg|\sum_{\{u,v\}\subset\Lambda}\bigg(\frac1{|\Lambda|}\bigg)^2
 \frac{\tau_{p_1}(u-v)}{g_{p_1}}-\frac12S_{\ge3}(o)\bigg|=O(\beta^2),
 \lbeq{H1''main1LA}\\
\bigg|\sum_{\{u,v\}\subset\Lambda}\bigg(\frac1{|\Lambda|}\bigg)^2\sum_{x\ne u,v}
 \frac{\tau_{p_1}(u-x)}{g_{p_1}}\frac{\tau_{p_1}(x-v)}{g_{p_1}}-\frac12S_{\ge
 2}^{*2}(o)\bigg|=O(\beta^2).\lbeq{H1''main2LA}
\end{gather}
\end{lmm}
\end{shaded}
We will prove Lemma~\ref{lmm:RWapproxLA} after the proof of Lemma~\ref{lmm:gp1a} is completed. 

Now we resume the proof of Lemma~\ref{lmm:gp1a} assuming the bounds in Lemma~\ref{lmm:RWapproxLA}.  First we recall 
$G=G_0+G_1+G_2$ (cf., \refeq{G-G0=G1+G2}), where $G_0$ is independent of the 
models and estimated as \refeq{G0-est}; $G_1$ is defined as \refeq{G1def} and 
here we use \refeq{RWest1LA} to show \refeq{G1-est}; $G_2$ is defined as 
\refeq{G2def} and obeys the same bound as \refeq{G2-est}.  As a result, 
Lemma~\ref{lmm:G} also holds for lattice animals.  Similarly, we can show 
$H=e(H_1'''+H_2'''-H_2''-H_2)+O(\beta^2)$ (cf., \refeq{H2bd}, 
\refeq{eH1'}, \refeq{H2''bd} and \refeq{H2'''bd}), 
where $H_1'''$ is defined in \refeq{H1'''def} and here we use 
\refeq{H1''main1LA}--\refeq{H1''main2LA} to show \refeq{H_1'''bd}; 
$H_2$ is bounded by $H_{2,3}+H_{2,4}$, and $H_{2,3}$ and $H_{2,4}$ are further 
bounded as \refeq{H2bd-pr2}--\refeq{H2bd-pr3} and 
\refeq{H2bd-pr4}--\refeq{H2bd-pr5} (with $\Tcal$ replaced by $\Acal$), 
and here we use the infrared bound \refeq{MFbehavior} and the convolution 
bound on power functions \cite[Lemma~3.2(i)]{cs15} to show $H_2=O(\beta^2)$, 
such as
\begin{align}
&\sum_{\{u,v\}\subset\Lambda}\bigg(\frac1{|\Lambda|}\bigg)^2\sum_{x\in\Zd}
 \tau_{p_1}(u-x)\,\tau_{p_1}(v-x)\nn\\
&\le\sum_{\{u,v\}\subset\Lambda}\bigg(\frac1{|\Lambda|}\bigg)^2\sum_{x\in\Zd}
 \frac{O(L^{-2})}{(\|u-x\|\vee L)^{d-2}}\frac{O(L^{-2})}{(\|v-x\|\vee L)^{d-2}}
 \nn\\
&\le\sum_{\{u,v\}\subset\Lambda}\bigg(\frac1{|\Lambda|}\bigg)^2\frac{O(L^{-4})}
 {(\|u-v\|\vee L)^{d-4}}~=O(\beta^2).
\end{align}
Similarly we can show that $H_2''$ and $H_2'''$ are both $O(\beta^2)$ by 
using the infrared bound and the convolution bound, instead of bounding 
$\tau_{p_1}/g_{p_1}$ by $S_{\ge0}$ or $S_{\ge1}$, just as done for lattice 
trees.
As a result, Lemma~\ref{lmm:H} also holds for lattice animals.

\begin{figure}[t]
\begin{center}
\includegraphics[scale=0.6]{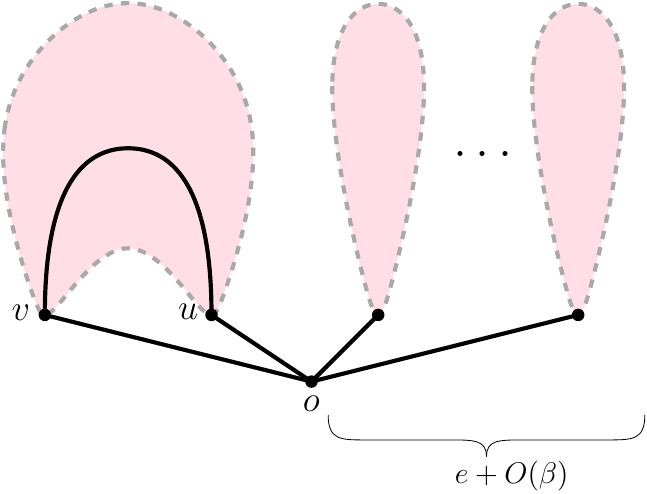}
\end{center}
\caption{Schematic representation of the dominant contribution to $I_1'$.}
\label{fig:I1p}
\end{figure}

\begin{figure}[t]
\[ \raisebox{-4pc}{\includegraphics[scale=0.5]{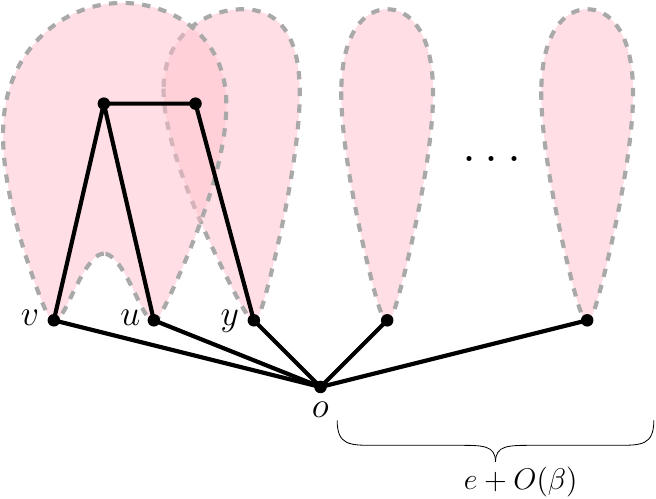}}\quad~+\qquad
 \raisebox{-4pc}{\includegraphics[scale=0.5]{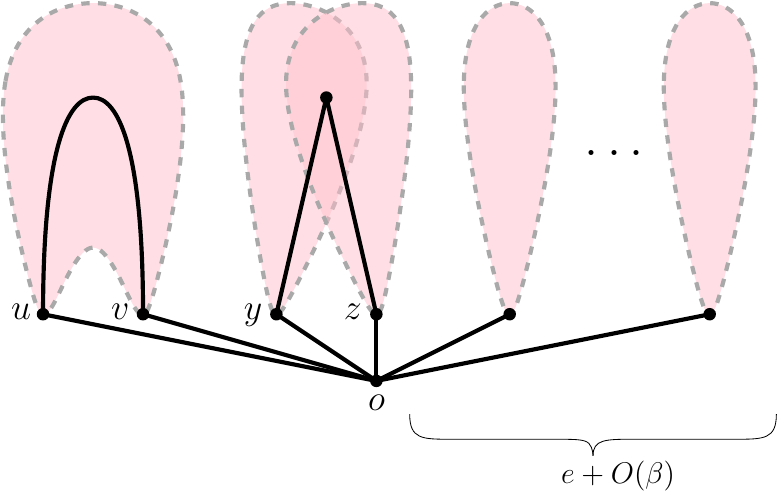}} \]
\caption{Schematic representation of the bound on \refeq{I1p-I1}.}
\label{fig:I234}
\end{figure}

Next we investigate $I$, which is unique for lattice animals.  
Let $I_1$ be the contribution from 
$\gamma\in\Gamma(Y)$ with $|\gamma|=|Y|-1$, i.e., consisting of a pair 
$\{u,v\}$ and $|Y|-2$ singletons $\{y\}_{y\in Y\setminus\{u,v\}}$:
\begin{gather}
I_1=\sum_{\substack{Y\subset\Lambda\\ (|Y|\ge2)}}\bigg(\frac{p_1}{|\Lambda|}
 \bigg)^{|Y|}\sum_{\{u,v\}\subset Y}\sum_{R\in\Acal_{u,v}\setminus\Acal_o}
 \bigg(\frac{p_1}{|\Lambda|}\bigg)^{|E_R|}\prod_{y\in Y\setminus\{u,v\}}
 \sum_{R_y\in\Acal_y\setminus\Acal_o}\bigg(\frac{p_1}{|\Lambda|}
 \bigg)^{|E_{R_y}|}\nn\\
\times\prod_{y\in Y\setminus\{u,v\}}\ind{V_R\cap V_{R_y}=\vno}
 \prod_{\substack{y,z\in Y\setminus\{u,v\}\\ (y\ne z)}}\ind{V_{R_y}\cap V_{R_z}
 =\vno},
\end{gather}
where an empty product is regarded as 1.  The dominant contribution to $I_1$, 
denoted $I'_1$, comes from when the last line is replaced by 1.  By the 
tree-graph inequality \refeq{bdof3pf}, which is also true for lattice animals due 
to subadditivity, and then using the infrared bound \refeq{MFbehavior}, 
it is estimated as (see Figure~\ref{fig:I1p})
\begin{align}\lbeq{I'_1}
I_1'&~=\sum_{\{u,v\}\subset\Lambda}\bigg(\frac{p_1}{|\Lambda|}\bigg)^2\Big(
 \tau_{p_1}(u-v)-\tau_{p_1}^{\sss(3)}(o,u,v)\Big)\underbrace{\sum_{Y'\subset
 \Lambda\setminus\{u,v\}}\bigg(\frac1{|\Lambda|}\bigg)^{|Y'|}\prod_{y\in Y'}
 \bigg(1-\frac{\tau_{p_1}(y)}{g_{p_1}}\bigg)}_{e+O(\beta)}\nn\\
&~=~\underbrace{p_1\big(e+O(\beta)\big)}_{1+O(\beta)}\sum_{\{u,v\}\subset
 \Lambda}\bigg(\frac1{|\Lambda|}\bigg)^2\bigg(\frac{\tau_{p_1}(u-v)}{g_{p_1}}
 -\frac{\tau_{p_1}^{\sss(3)}(o,u,v)}{g_{p_1}}\bigg)\nn\\
&~=\sum_{\{u,v\}\subset\Lambda}\bigg(\frac1
 {|\Lambda|}\bigg)^2\frac{\tau_{p_1}(u-v)}{g_{p_1}}+O(\beta^2)\nn\\
&\stackrel{\text{\refeq{H1''main1LA}}}=~\frac12S_{\ge3}(o)+O(\beta^2).
\end{align}
On the other hand, by using $1-ab\le(1-a)+(1-b)$ for any $a,b\in\{0,1\}$, 
we can bound the difference $I'_1-I_1~(\ge0)$ as
\begin{align}\lbeq{I1p-I1}
I'_1-I_1
&\le\sum_{\{u,v\}\subset\Lambda}\bigg(\frac{p_1}{|\Lambda|}\bigg)^2\sum_{R\in
 \Acal_{u,v}\setminus\Acal_o}\bigg(\frac{p_1}{|\Lambda|}\bigg)^{|E_R|}\!\!
 \sum_{Y'\subset\Lambda\setminus\{u,v\}}\bigg(\frac{p_1}{|\Lambda|}\bigg)^{|Y'|}
 \prod_{y\in Y'}\sum_{R_y\in\Acal_y\setminus\Acal_o}\bigg(\frac{p_1}{|\Lambda|}
 \bigg)^{|E_{R_y}|}\nn\\
&\hskip5pc\times\bigg(1-\prod_{y\in Y'}\ind{V_R\cap V_{R_y}=\vno}+1-\prod_{
 \substack{y,z\in Y'\\ (y\ne z)}}\ind{V_{R_y}\cap V_{R_z}=\vno}\bigg)\nn\\
&\le\sum_{\{u,v\}\subset\Lambda}\bigg(\frac{p_1}{|\Lambda|}\bigg)^2\sum_{R\in
 \Acal_{u,v}\setminus\Acal_o}\bigg(\frac{p_1}{|\Lambda|}\bigg)^{|E_R|}\!\!
 \sum_{Y'\subset\Lambda\setminus\{u,v\}}\bigg(\frac{p_1}{|\Lambda|}\bigg)^{|Y'|}
 \prod_{y\in Y'}\sum_{R_y\in\Acal_y\setminus\Acal_o}\bigg(\frac{p_1}{|\Lambda|}
 \bigg)^{|E_{R_y}|}\nn\\
&\hskip5pc\times\bigg(\sum_{y\in Y'}\ind{V_R\cap V_{R_y}\ne\vno}+\sum_{
 \substack{y,z\in Y'\\ (y\ne z)}}\ind{V_{R_y}\cap V_{R_z}\ne\vno}\bigg).
\end{align}
This is $O(\beta^2)$, as the contribution from the former (resp., latter) 
sum in the last line can be estimated in a similar way to showing 
$H_{2,3}=O(\beta^2)$ (resp., $H_{2,4}=O(\beta^2)$); see Figure~\ref{fig:I234}.  
As a result, 
\begin{align}
I_1=\frac12S_{\ge3}(o)+O(\beta^2).
\end{align}

Finally we estimate the difference $I-I_1$:
\begin{align}\lbeq{I-I1}
I-I_1&=\sum_{\substack{Y\subset\Lambda\\ (|Y|\ge3)}}\bigg(\frac{p_1}{|\Lambda|}
 \bigg)^{|Y|}\sum_{\substack{\gamma\in\Gamma(Y)\\ (|\gamma|\le|Y|-2)}}\prod_{j
 =1}^{|\gamma|}\sum_{R_j\in\Acal_{\gamma_j}\setminus\Acal_o}\bigg(\frac{p_1}
 {|\Lambda|}\bigg)^{|E_{R_j}|}\prod_{i<j}\ind{V_{R_i}\cap V_{R_j}=\vno}.
\end{align}
Since $|\gamma|\le|Y|-2$, there are two possibilities: (i) there is a set in 
$\gamma$ which includes at least 3 distinct neighbors of the origin, or 
(ii) there are at least two disjoint sets in $\gamma$ both of which include 
exactly two distinct neighbors of the origin.  Therefore,
\begin{align}
I-I_1&\le\sum_{\substack{Y\subset\Lambda\\ (|Y|\ge3)}}\bigg(\frac{p_1}
 {|\Lambda|}\bigg)^{|Y|}\sum_{\gamma\in\Gamma(Y)}\prod_{j=1}^{|\gamma|}
 \sum_{R_j\in\Acal_{\gamma_j}\setminus\Acal_o}\bigg(\frac{p_1}{|\Lambda|}
 \bigg)^{|E_{R_j}|}\prod_{i<j}\ind{V_{R_i}\cap V_{R_j}=\vno}\nn\\
&\hskip5pc\times\Big(\ind{\exists j,\;|\gamma_j|\ge3}+\ind{\exists i\ne j,\;
 |\gamma_i|=|\gamma_j|=2}\Big)\nn\\
&=I_3+I_2,
\end{align}
where $I_3$ and $I_2$ are the contributions from 
$\ind{\exists j,\;|\gamma_j|\ge3}$ and 
$\ind{\exists i\ne j,\;|\gamma_i|=|\gamma_j|=2}$, respectively.

For $I_2$, we split the set $Y$ of neighbors of the origin into $U$, $V$ and 
$Y'=Y\setminus(U\cup V)$, where $U\cap V=\vno$ and $|U|=|V|=2$.  
Partially ignoring the avoidance constraint among animals, we can bound $I_2$ as
\begin{align}
I_2&\le\sum_{\substack{U\subset\Lambda \\ (|U|=2)}}\bigg(\frac{p_1}{|\Lambda|}
 \bigg)^2\sum_{R\in\Acal_U\setminus\Acal_o}\bigg(\frac{p_1}{|\Lambda|}
 \bigg)^{|E_R|}\sum_{\substack{V\subset\Lambda\setminus U \\ (|V|=2)}}\bigg(
 \frac{p_1}{|\Lambda|}\bigg)^2\sum_{R'\in\Acal_V\setminus\Acal_o}\bigg(
 \frac{p_1}{|\Lambda|}\bigg)^{|E_{R'}|}\nn\\
&\qquad\times\sum_{Y'\subset\Lambda\setminus(U\cup V)}\bigg(\frac{p_1}
 {|\Lambda|}\bigg)^{|Y'|}\sum_{\gamma\in\Gamma(Y')}\prod_{j=1}^{|\gamma|}
 \sum_{R_j\in\Acal_{\gamma_j}\setminus\Acal_o}\bigg(\frac{p_1}{|\Lambda|}
 \bigg)^{|E_{R_j}|}\prod_{i<j}\ind{V_{R_i}\cap V_{R_j}=\vno}.
\end{align}
Notice that the second line is almost identical to $g_{p_1}$; the only 
difference is the domain of summation over $Y'$, and therefore it is bounded 
above by $g_{p_1}$.  Since $p_1g_{p_1}=1$, we obtain
\begin{align}
I_2\le g_{p_1}\Bigg(\sum_{\{u,v\}\subset\Lambda}\bigg(\frac{p_1}{|\Lambda|}
 \bigg)^2\underbrace{\sum_{R\in\Acal_{u,v}}\bigg(\frac{p_1}{|\Lambda|}
 \bigg)^{|E_R|}}_{\tau_{p_1}(u-v)}\Bigg)^2
 \stackrel{\text{\refeq{H1''main1LA}}}\le p_1\bigg(\frac12S_{\ge3}(o)
 +O(\beta^2)\bigg)^2=O(\beta^2).
\end{align}

For $I_3$, we split the set $Y$ into $X$ and $Y'=Y\setminus X$, where $X$ 
includes at least 3 distinct vertices $x,y,z\in\Lambda$.  Again, by partially 
ignoring the avoidance constraint among animals, we can bound $I_3$ as
\begin{align}
I_3&\le\sum_{\{x,y,z\}\subset\Lambda}\sum_{\substack{X
 \subset\Lambda \\ (X\ni x,y,z)}}\bigg(\frac{p_1}{|\Lambda|}\bigg)^{|X|}\sum_{R
 \in\Acal_X\setminus\Acal_o}\bigg(\frac{p_1}{|\Lambda|}\bigg)^{|E_R|}\nn\\
&\qquad\times\sum_{Y'\subset\Lambda\setminus X}\bigg(\frac{p_1}
 {|\Lambda|}\bigg)^{|Y'|}\sum_{\gamma\in\Gamma(Y')}\prod_{j=1}^{|\gamma|}
 \sum_{R_j\in\Acal_{\gamma_j}\setminus\Acal_o}\bigg(\frac{p_1}{|\Lambda|}
 \bigg)^{|E_{R_j}|}\prod_{i<j}\ind{V_{R_i}\cap V_{R_j}=\vno}.
\end{align}
Notice again that the second line is bounded above by $g_{p_1}$.  
Using the relation 
$\Acal_X\setminus\Acal_o\subset\Acal_{x,y,z}$ and splitting $X$ into 
$\{x,y,z\}$ and $X'=X\setminus\{x,y,z\}$, we obtain
\begin{align}
I_3&\le g_{p_1}\sum_{\{x,y,z\}\subset\Lambda}\bigg(
 \frac{p_1}{|\Lambda|}\bigg)^3\underbrace{\sum_{R\in\Acal_{x,y,z}}\bigg(
 \frac{p_1}{|\Lambda|}\bigg)^{|E_R|}}_{\tau_{p_1}^{(3)}(x,y,z)}~\underbrace{
 \sum_{X'\subset\Lambda\setminus\{x,y,z\}}\bigg(\frac{p_1}{|\Lambda|}
 \bigg)^{|X'|}}_{\le\,(1+p_1/|\Lambda|)^{|\Lambda|}}.
\end{align}
By the tree-graph inequality \refeq{bdof3pf}, we can show that
\begin{align}
&\sum_{\{x,y,z\}\subset\Lambda}\bigg(\frac{p_1}{|\Lambda|}\bigg)^3
 \tau_{p_1}^{\sss(3)}(x,y,z)\nn\\
&\le\sum_{\{x,y,z\}\subset\Lambda}\bigg(\frac1{|\Lambda|}\bigg)^3\sum_{w\in\Zd}
 \frac{\tau_{p_1}(x-w)}{g_{p_1}}\frac{\tau_{p_1}(y-w)}{g_{p_1}}\frac{\tau_{p_1}
 (z-w)}{g_{p_1}}\nn\\
&=\sum_{\{x,y,z\}\subset\Lambda}\bigg(\frac1{|\Lambda|}\bigg)^3\bigg(\sum_{w
 \ne x,y,z}\frac{\tau_{p_1}(x-w)}{g_{p_1}}\frac{\tau_{p_1}(y-w)}{g_{p_1}}
 \frac{\tau_{p_1}(z-w)}{g_{p_1}}+\frac3{g_{p_1}}\frac{\tau_{p_1}(x-z)}{g_{p_1}}
 \frac{\tau_{p_1}(y-z)}{g_{p_1}}\bigg)\nn\\
&\le\frac3{g_{p_1}}\Big(\underbrace{\|D*\tau_{p_1}\|_\infty+\|D\|_\infty}_{O
 (\beta)~(\because\text{\refeq{MFbehavior}})}\Big)\underbrace{\sum_{\{x,y\}
 \subset\Lambda}\bigg(\frac1{|\Lambda|}\bigg)^2\sum_{w\ne x,y}\frac{\tau_{p_1}
 (x-w)}{g_{p_1}}\frac{\tau_{p_1}(y-w)}{g_{p_1}}}_{O(\beta)~(\because
 \text{\refeq{H1''main2LA}}},
\end{align}
hence $I_3=O(\beta^2)$.  This completes the proof of 
$I=\frac12S_{\ge3}(o)+O(\beta^2)$, hence the proof of Lemma~\ref{lmm:gp1a}.
\QED

\Proof{Proof of Lemma~\ref{lmm:RWapproxLA}.}
First we prove \refeq{RWest1LA}.  By the inverse Fourier transfom, we have 
the rewrite
\begin{align}\lbeq{convoDtau}
\sum_{u\in\Lambda}\frac1{|\Lambda|}\frac{\tau_{p_1}(u)}{g_{p_1}}=\bigg(D*
 \frac{\tau_{p_1}(u)}{g_{p_1}}\bigg)(o)=\int_{[-\pi,\pi]^d}\hat D(k)\,\frac{\hat
 \tau_{p_1}(k)}{g_{p_1}}\frac{\diff^dk}{(2\pi)^d}.
\end{align}
Notice that the Fourier transform of the recursion equation \refeq{lace-exp} yields 
\begin{align}
\hat \tau_p(k)
 =\frac{g_p+\hat h_p(k)+\hat\pi_p(k)}{1-(g_p+\hat h_p(k)+\hat\pi_p(k))p\hat
 D(k)}.
\end{align}
We use this identity at $p_1=1/g_{p_1}$.  Let 
\begin{align}
H(x)=\frac{h_{p_1}(x)+\pi_{p_1}(x)}{g_{p_1}}.
\end{align}
Thanks to the symmetry, the Fourier transform $\hat H(k)$ is real.  
Moreover, by \refeq{pin-hatbd}--\refeq{hathbd}, we can show that, 
for $d>8$ and $L\gg1$, $|\hat H(k)|=O(\beta)$ uniformly in $k$.  
Then, we can rewrite $\hat\tau_{p_1}(k)/g_{p_1}$ as
\begin{align}\lbeq{hatFdef}
\frac{\hat\tau_{p_1}(k)}{g_{p_1}}
&=\frac{1+\hat H(k)}{1-(1+\hat H(k))\hat D(k)}\nn\\
&=\frac1{1-\hat D(k)}+\frac{\hat H(k)}{1-\hat D(k)}\underbrace{\frac1{1-(1+\hat
 H(k))\hat D(k)}}_{=:\hat F(k)}.
\end{align}
Applying this to \refeq{convoDtau} yields the main term $S_{\ge2}(o)$ as
\begin{align}\lbeq{FTD*tau}
\sum_{u\in\Lambda}\frac1{|\Lambda|}\frac{\tau_{p_1}(u)}{g_{p_1}}
&=\underbrace{\int_{[-\pi,\pi]^d}\frac{\hat D(k)}{1-\hat D(k)}\frac{\diff^dk}
 {(2\pi)^d}}_{\hskip4pc S_{\ge1}(o)~(=S_{\ge2}(o))}+\int_{[-\pi,\pi]^d}
 \frac{\hat D(k)\hat H(k)}{1-\hat D(k)}\hat F(k)\frac{\diff^dk}{(2\pi)^d}.
\end{align}

It remains to show that the second term on the right is $O(\beta^2)$.  
To do so, we want an effective bound on $\hat F(k)$.  We will show at the end 
of the proof that, for $d>8$ and $L\gg1$, there is an $L$-independent constant 
$C<\infty$ such that
\begin{align}\lbeq{hatFnaivebd}
0<\hat F(k)\le\frac{C}{1-\hat D(k)},
\end{align}
uniformly in $k$. However, to use \refeq{hatFnaivebd} for the second term of \refeq{FTD*tau}, we have to bound $\hat D(k)$ above by the absolute value of it, which makes it difficult to show the error being $O(\beta^2)$. Instead, 
we first rewrite $\hat F(k)$ as
\begin{align}\lbeq{FTtaup1-2}
\hat F(k)=\frac1{1-\hat D(k)}+\frac{\hat D(k)\hat H(k)}{1-\hat D(k)}\hat F(k).
\end{align}
Then, the second term on the right of \refeq{FTD*tau} equals
\begin{align}\lbeq{DSHF}
\underbrace{\int_{[-\pi,\pi]^d}\frac{\hat D(k)\hat H(k)}{(1-\hat D(k))^2}
 \frac{\diff^dk}{(2\pi)^d}}_{(D*S_{\ge0}^{*2}*H)(o)}+\int_{[-\pi,\pi]^d}\bigg(
 \frac{\hat D(k)\hat H(k)}{1-\hat D(k)}\bigg)^2\hat F(k)\frac{\diff^dk}
 {(2\pi)^d}.
\end{align}
Notice that, due to the identity \refeq{FTtaup1-2}, we can make $\hat D(k)^2$, which is always positive, in the second term of \refeq{DSHF}. The first term  is readily bounded by 
$\|D*S_{\ge0}^{*2}\|_\infty\hat H(0)=O(\beta^2)$.  For the second 
term, we use $|\hat H(k)|=O(\beta)$ and \refeq{hatFnaivebd} 
to obtain that
\begin{align}\lbeq{naivebd-appl}
\int_{[-\pi,\pi]^d}\bigg(\frac{\hat D(k)\hat H(k)}{1-\hat D(k)}\bigg)^2\hat F(k)
 \frac{\diff^dk}{(2\pi)^d}&\le O(\beta^2)\underbrace{\int_{[-\pi,\pi]^d}
 \frac{\hat D(k)^2}{(1-\hat D(k))^3}\frac{\diff^dk}{(2\pi)^d}}_{(D^{*2}*S_{\ge
 0}^{*3})(o)}=O(\beta^3).
\end{align}
This completes the proof of \refeq{RWest1LA}.

We can also prove \refeq{H1''main1LA}--\refeq{H1''main2LA} in a similar manner by assuming \refeq{hatFnaivebd}. Hence we here prove only \refeq{H1''main1LA}.
By the inverse Fourier transform, we can rewrite the sum in 
\refeq{H1''main1LA} as
\begin{align}
\sum_{\{u,v\}\subset\Lambda}\bigg(\frac1{|\Lambda|}\bigg)^2\frac{\tau_{p_1}
 (u-v)}{g_{p_1}}&=\frac12\left(D^{*2}*\frac{\tau_{p_1}}{g_{p_1}} \right)(o)
 -\frac12D^{*2}(o)\nn\\
&=\frac12\int_{[-\pi,\pi]^d}\hat D(k)^2\frac{\hat\tau_{p_1}(k)}{g_{p_1}}
 \frac{\diff^dk}{(2\pi)^d}-\frac12D^{*2}(o).
\end{align}
Then, by the identity \refeq{hatFdef}, we can extract the main term 
$\frac12S_{\ge3}(o)$ as
\begin{align}
\underbrace{\frac12\int_{[-\pi,\pi]^d}\frac{\hat D(k)^2}{1-\hat D(k)}
 \frac{\diff^dk}{(2\pi)^d}-\frac12D^{*2}(o)}_{\frac12S_{\ge3}(o)}+\,\frac12
 \int_{[-\pi,\pi]^d}\frac{\hat D(k)^2\hat H(k)}{1-\hat D(k)}\hat F(k)
 \frac{\diff^dk}{(2\pi)^d}.
\end{align}
Similarly to \refeq{naivebd-appl}, the second term is bounded as 
\begin{align}
\bigg|\frac12\int_{[-\pi,\pi]^d}\frac{\hat D(k)^2\hat H(k)}{1-\hat D(k)}\hat
 F(k)\frac{\diff^dk}{(2\pi)^d}\bigg|\stackrel{\text{\refeq{hatFnaivebd}}}\le
 O(\beta)\underbrace{\int_{[-\pi,\pi]^d}\frac{\hat D(k)^2}{(1-\hat D(k))^2}
 \frac{\diff^dk}{(2\pi)^d}}_{(D^{*2}*S_{\ge0}^{*2})(o)}=O(\beta^2),
\end{align}
hence the completion of the proof of \refeq{H1''main1LA}.

Finally we prove the inequality \refeq{hatFnaivebd}, for $\|k\|\ge\frac1L$ and 
$\|k\|\le\frac1L$ separately.  We begin with the former case.  It is known (cf., 
e.g., \cite{hs05}) that our $D$ satisfies \cite[Assumption~D]{hs02}; in 
particular, there is an $L$-independent constant $\eta\in(0,1)$ such that
\begin{align}\lbeq{k>1/L}
-1+\eta~\stackrel{\forall k}\le~\hat D(k)\stackrel{\|k\|\ge\frac1L}\le1-\eta.
\end{align}
Since $|\hat H(k)|=O(\beta)$, we obtain that, for $L\gg1$,
\begin{align}
-1+\frac\eta2~\stackrel{\forall k}\le~\big(1+\hat H(k)\big)\hat D(k)
 \stackrel{\|k\|\ge\frac1L}\le1-\frac\eta2,
\end{align}
hence
\begin{align}\lbeq{hatFbd4largek}
0<\frac1{2-\eta/2}~\stackrel{\forall k}\le~\hat F(k)\stackrel{\|k\|\ge\frac1L}
 \le\frac2\eta~\stackrel{\forall k}\le~\frac2\eta\frac{2-\eta}{1-\hat D(k)}.
\end{align}

It remains to show that $\hat F(k)$ is bounded above by a multiple of 
$(1-\hat D(k))^{-1}$ uniformly in $\|k\|\le\frac1L$.  We note that
\begin{align}\lbeq{1/hatFrewr}
\hat F(k)^{-1}&=1-\big(1+\hat H(0)\big)+\big(1+\hat H(0)\big)\big(1-\hat D(k)
 \big)+\big(\hat H(0)-\hat H(k)\big)\hat D(k)\nn\\
&=-\hat H(0)+\bigg(1+\hat H(0)+\frac{\hat H(0)-\hat H(k)}{1-\hat D(k)}\hat
 D(k)\bigg)\big(1-\hat D(k)\big).
\end{align}
Since $-\hat H(0)$ is bounded below by  a positive multiple of $\beta$ 
(as explained in the beginning of the proof of Lemma~\ref{lmm:pc-2nd}), 
ignoring this term yields a lower bound on $\hat F(k)^{-1}$.  Moreover, since 
$|k\cdot x|\le\|k\|\|x\|\le1$ for $x\in\Lambda$ 
and $\|k\|\le\frac1L$, and since $1-\cos t\ge\frac2{\pi^2}t^2$ for $|t|\le1$, 
there is a $c>0$ such that
\begin{align}
1-\hat D(k)=\sum_{x\in\Lambda}\frac{1-\cos(k\cdot x)}{|\Lambda|}\ge\frac2
 {\pi^2}\sum_{x\in\Lambda}\frac{(k\cdot x)^2}{|\Lambda|}=\frac{2\|k\|^2}{d\pi^2}
 \underbrace{\sum_{x\in\Lambda}\frac{\|x\|^2}{|\Lambda|}}_{\ge cL^2}.
\end{align}
On the other hand, by $1-\cos t\le\frac12t^2$ for any $t$ and using the 
$x$-space bounds \refeq{pin-xbd} and \refeq{hathbd}, we have 
\begin{align}
|\hat H(0)-\hat H(k)|\le\sum_x\frac{(k\cdot x)^2}2|H(x)|\le\frac{\|k\|^2}{2d}
 \underbrace{\sum_x\|x\|^2|H(x)|}_{O(L^2\beta)}.
\end{align}
Therefore, by taking $L$ sufficiently large, $\hat F(k)^{-1}$ is bounded 
below by a positive multiple of  $1-\hat D(k)$, uniformly in $\|k\|\le\frac1L$. 
Combined with \refeq{hatFbd4largek}, this completes the proof of the inequality 
\refeq{hatFnaivebd}, hence the completion of the proof of 
Lemma~\ref{lmm:RWapproxLA}.
\QED

\section*{Acknowledgements}
This work was supported by JSPS KAKENHI Grant Number 18K03406.  We are 
grateful to Yinshuang Liang for working together in the early stage of this 
project.  We would also like to thank Gordon Slade for comments to an earlier 
version of this paper. We are grateful to the two anonymous referees for thorough reviewing and many valuable suggestions to the previous version of the manuscript.

\end{document}